\documentclass[12pt]{amsart}
\usepackage {amsmath, amscd, amsbsy, amsfonts}
\usepackage{amsmath,latexsym,amssymb,amsfonts,mathrsfs,bbm}
\usepackage{fullpage}
\usepackage{graphicx}
\usepackage{hyperref}
\newlength{\baseunit}               
\newcommand{\R}{\mathbbmss{R}}

\def \setA {\{A\}}
\def \Y {{\mathcal Y}}
\def \gam {{\Gamma}}
\def \seta {\{A\}}
\def \n {{\mathcal N}}

\def \D {\mathcal D}
\def \Q {\mathbb Q}
\def \proj {{\mathbb{P}^r}}

\def \X {\mathcal X}

\def \d {\Delta}

\def \rationals {\mathbb{Q}}
\def \L {{\mathcal L}}
\def \W {{\mathcal W}}

\def \H {{\mathcal H}}
\def \K {{\mathcal K}}
\def \V {{\mathcal V}}
\def \P {{\mathcal P}}
\def \G {{\mathcal G}}

\def \F {{\mathcal F}}

\def \R {{\mathcal R}}
\def \C {{\mathcal C}}

\def \T {{\mathcal T}}
\def \n {{\mathcal N}}
\def \mbar {\overline{\mathcal M}}
\def \sep {\ || \ }
\def \wt {\widetilde}

\def \V {{\mathcal V}}

\def \J {{\mathcal J}}

\newtheorem{theorem}{Theorem}[section] 

\newtheorem{lemma}[theorem]{Lemma}
\newtheorem{example}[theorem]{Example}
\newtheorem{corollary}[theorem]{Corollary}

\newtheorem{proposition}[theorem]{Proposition} 
 
\newcommand{\bpf}{ {\it Proof.  }}
\newcommand{\epf}{\qed \vspace{+10pt}}

\title[Enumerative Invariants.]{Characteristic numbers of elliptic curves with fixed j-invariant}

\author{Dung Nguyen}

\begin{document}
\maketitle
\begin{abstract}
We solve the problem of counting elliptic curves with 
fixed j-invariant in projective space with tangency conditions.
This is equivalent to couting rational nodal curves with
condition on the node of the image. The solution is given 
in the form of effective recursions. We give explicit formulas when the dimension
of the ambient projective space is at most $5$. Many
numerical examples are provided. A C++ program
implementing all of the recursions 
is available upon request.
 \vspace{-20pt}
\end{abstract} 
\section{introduction}
 Charateristic numbers of curves in projective spaces is a classical problem in algebraic
geometry: how many curves in $\proj$ of given degree and genus 
that pass through a general set
of linear subspaces, and are tangent to a general set of hyperplanes? 
Presented in this form, the problem seems almost
unattackable, as not much is known even in the case of genus two space curves.
However, the cases of genus zero and genus one space curves are well
understood. Incidence-only (meaning no tangency condition is considered)
characteristic numbers of rational plane curves
 were first computed by Kontsevich, see \cite{fp}. The method was to pull back the WDVV 
equation on $\mbar_{0,4}$ onto the moduli space of stable maps
$\mbar_{0,n}(2,d)$ to obtain a recursion counting rational plane
curves. The same method works equally well for rational space curves. In \cite{idq}, 
Lemma $2.3.1$,
it was shown that the tangency divisor is numerically equivalent
to a linear combination of the incident divisor and boundary
divisors on $\mbar_{0,n}(r,d).$ Hence one can write down
a recursion computing full characteristic numbers of
rational space curves. 

  In genus one, there are at least two counting problems.
 One could try to obtain enumeration of genus one curves with
generic $j-$ invariant, or of genus one curves with fixed $j-$
invariant. This note will deal with the latter. 
Incidence-only characteristic numbers
for genus one space curves with fixed-j invariant have
been computed in \cite{eln} and \cite{zin}. In this
note, recursions computing all characteristic numbers will be 
provided. In case of incidence-only numbers, we obtain an algebraic
solution that works over any closed field of zero charactersistic,
in contrast to the analytic method in \cite{eln} and
\cite{zin}. The results in this note will also be used
to compute characteristic numbers of elliptic space curves
in an upcoming paper by the author.

All the recursions are based on our algorithm counting rational 
two nodal reducible curves. These are projective curves having 
two rational smooth component intersecting at two points (or with
a choice of two intersection points in the case of plane curves). 
Counting these curves is in turn based on an algorithm counting rational curves,
 now with an additional type of conditions: special tangent conditions.
 This will be defined in Section $2$.
 We work out in detail the algorithm counting rational curves with 
special tangent conditions in ambient space of dimension at most $5$.
 For dimension $6$ or higher, the numbers could in theory be expressed as
intersections of tautological classes on a blowup of
$\mbar_{0,1}(r,d)$, but this is much less implementable.

We use the following results to obtain our recursions. We use
the WDVV equation on $\mbar_{0,n}(r,d)$. We use the divisor
theory on $\mbar_{0,n}(r,d)$ as developed in \cite{idq}. We
do not use any outside input, and our method for incidence-only
characteristic numbers is different
from those in \cite{eln}, \cite{zin}. 

The author is very grateful to R.Vakil, his advisor, for numerous
helpful conversations and ideas, and for introducing him
to the beautiful subject of enumerative geometry.
\\ \\

\section{Definitions and Notations}
\subsection{The moduli space of stable maps of genus $0$ in $\proj$.}
As usual, $\mbar_{0,n}(r,d)$ will denote the Kontsevich compactification of the moduli
space of genus zero curves with $n$ marked points of degree $d$ in $\proj$. 
 We will also be using the notation $\mbar_{0,S}(r,d)$ where the markings 
are indexed by a set $S$.
The following are Weil divisors on $\mbar_{0,S}(r,d)$: 
\begin{itemize}
 \item The divisor $(U \sep V)$ of $\mbar_{0,S}(r,d)$
 is the closure in $\mbar_{0,S}(r,d)$ of the locus 
of curves with two components such that $U \cup V = S$ is a partition 
of the marked points over the two components.
\item The divisor $(d_1,d_2)$ is the closure in $\mbar_{0,S}(r,d)$
 of the locus of curves with two components, sucht that $d_1 + d_2 = d$ 
is the degree partition over the two components. 
\item The divisor $(U,d_1 \sep V,d_2)$ 
is the closure in $\mbar_{0,S}(r,d)$ of the locus of 
curves with two components, where $U \cup V = S$
and $d_1 +d_2 = d$ are the partitions of markings 
and degree over the two components respectively.
\end{itemize}

\subsection{The constraints and the ordering of constraints.}

 We will be concerned with the number of curves passing through a constraint. Each 
constraint is denoted by a $(r+1)-$tuple $\d$ as follows : \\
{\bf (i)}  $\d(0)$ is the number of hyperplanes that the curves need to be tangent to. \\
{\bf (ii)} For $0<i \leq r$, $\d(i)$ is the number of subspaces of codimension
$i$ that the curves need to pass through. \\
{\bf (iii)} If the curves in consideration have a node and we place a condition
on the node, that is the node has to belong to a general codimension
$k$ linear subspace, then $\d$ has $r+2$ elements and the last element, $\d(r+1)$, is $k$.

 Note that because in general a curve of degree
$d$ will always intersect a hyperplane at $d$ points, introducing an incident condition
with a hyperplane essentially means multiplying the cycle class cut out
by other conditions by $d$.
For example, if we ask how many genus zero curves of degree $4$ in $\mathbb P^3$ that pass through the constraint
$\d = (1,2,3,4,0) (\d(1) = 2)$, that means we ask how many genus zero curves of degree $4$ pass through
three lines, four points, are tangent to one hyperplane, and then multiply that answer
by $4^2$. We will also refer to $\d$ as a set of linear spaces,
hence we can say, pick a space $p$ in $\d$. 

 We consider the following ordering on the set of constraints, in order to prove that our algorithm will 
terminate later on. Let $r(\d) = -\sum_{i >1}^{i\leq r} \d[i]\cdot i^2$, and this will be our rank function.
 We compare two constraints $\d,\d'$ using the 
following criteria, whose priority are in the following order :
\begin{itemize}
 \item  If $\d(0) = \d'(0)$ and $\d$ has fewer non-hyperplane elements than $\d'$ does, then $\d<\d'$. 
 \item If $\d(0) > \d'(0)$  then $\d < \d'$.
\item If $r(\d) < r(\d')$ then $\d < \d'$. 
\end{itemize}

Informally speaking, characteristic numbers
where the constraints are more spread out at two ends 
are computed first in the recursion.
We write $\d = \d_1\d_2$ if $\d = \d_1 + \d_2$
as a parition of the set of linear spaces in $\d.$

\subsection{The stacks $\R,\n, \R\R, \R\R_2$ .}
 We list the following definitions of stacks of stable maps that will
occur in our recursions. \\ \\
{\bf 1)} Let $\R(r,d)$  be the usual moduli space of genus zero stable maps
$\mbar_{0,0}(r,d)$. \\ \\
{\bf 2)} Let $\n(r,d)$ be the closure in $\mbar_{0,\{A,B\} }(r,d)$ of the locus of maps of smooth
rational curves $\gamma$ such that $\gamma(A) = \gamma(B)$. Informally,
$\n(r,d)$ parametrizes degree $d$ rational nodal curves in $\proj$. \\ \\
{\bf 3)} For $d_1,d_2>0,$ let $\R\R(r,d_1,d_2)$ be $\mbar_{0,\{C\}}(r,d_1) \times \mbar_{0,\{C\}}(r,d_2)$
where the fibre product is taken over evaluation maps $ev_{C}$ to $\proj.$ \\ \\
{\bf 4)} Similarly we can define $\n\R(r,d_1,d_2)$ (see figure 1).\\ \\
{\bf 5)} For $d_1,d_2>0$, let $\R\R_2(r,d_1,d_2)$ be the closure in $\mbar_{0,\{A,C\}}(r,d_1) \times_{\proj} \mbar_{0,\{B,C\}}(r,d_2)$  
(the projections are evaluation maps $e_C$) of the locus of maps $\gamma$ 
such that $\gamma(A) = \gamma(B)$. We call maps in this family
rational two-nodal reducible curves.
\\ \\ 
$$\includegraphics[width = 150mm]{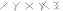}$$ \\
$$\text{Fig 1. Pictorial description of a general curve in the stacks $\R,\n, \R\R, \n\R, \R\R_2$}$$

\subsection{Special Tangent Condition} It is necessary to understand the enumerative geometry of
rational curves, now considering extra conditions of the form: there is a fixed
marked point $A$ on the curve, and the projective tangent line
at $A$ passes through a given codimension $2$ linear subspace $M$. 
The corresponding (Weil) divisor is denoted by $\W_A^M$. When
there is no need to consider any particular codimension $2$ subspace
$M$, we will only write $\W_A$.
We would also need to
consider the case where there is a condition on $A$, which means it could be
specified to lie on a certain linear subspace. By characteristic numbers
of rational space curves with special tangent conditions, we mean the numbers
of rational space curves having a marked point $A$ that satisfy the following conditions :
\begin{itemize}
\item Pass through various linear spaces and are tangent
to various hyperplanes.
\item The tangent line at $A$ to the curve passes through various codimension $2$ linear spaces.
\item The point $A$ may or may not lie on a given linear space.
\end{itemize}
$$\includegraphics[width = 60mm]{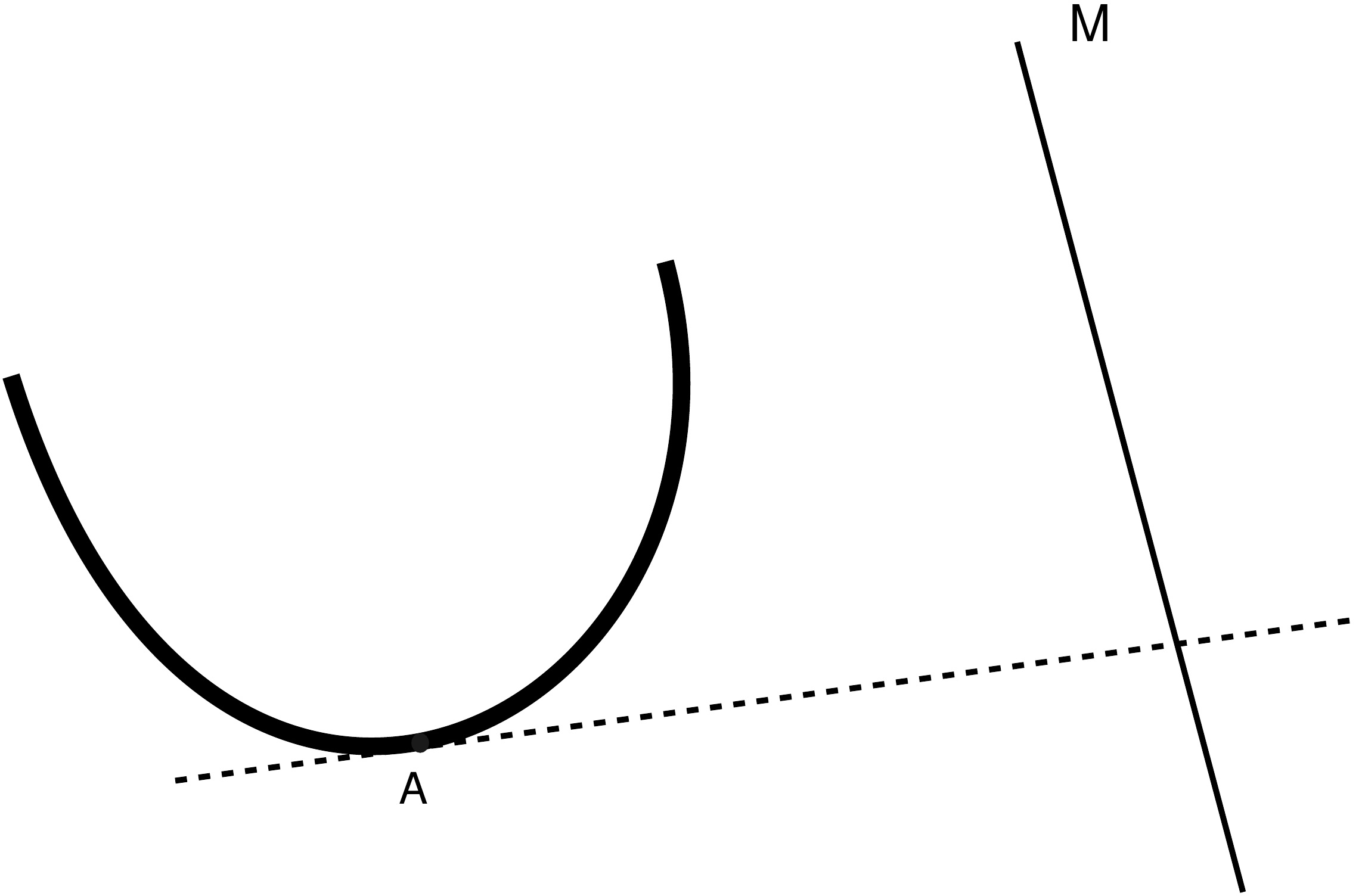}$$
$$\text{Fig 2. A curve with a special tangent condition}$$

\subsection{Stacks of stable maps with constraints.} Let $\F$
be a 
 maps of curves into $\proj$. For a constraint
$\d$, we define
$(\F, \d)$ be the closure in $\F $ of the locus of maps
that satisfy the constraint $\d$. If the stack of maps
$\F$ has two marked points $A$ and $B$, we define
$(\F, \L_A^u\L_B^v)$ to be the closure
in $\F$ of the locus of maps $\gamma$ such 
that $\gamma(A)$ lies on $u$ general hyperplanes,
 and that  $\gamma(B)$ lies on $v$ general hyperplanes. 

If $\F$ has one marked point $A$ then 
we define $(\F,\L_A^u\W^v_A)$ 
to be the closure of maps $\gamma$ such that
$\gamma(A)$ lies on $u$ general hyperplanes, and that
the image of $\gamma$ is smooth at $\gamma(A)$ and the tangent line to the image of $\gamma$
at $\gamma(A)$ passes through $v$ general codimension $2$
subspaces ($v$ special tangent conditions)

If a stack $\F$ is supported on a finite
number of points then we denote
$\# \F$ to be the stack-theoretic length of $\F$.

If $\F$ is a closed substack of the stacks $\n\R, \R\R$ then we denote
$(\F,\gam_1,\gam_2,k)$  to be the closure in $\F$
of the locus of maps $\gamma$ such that the
restriction of $\gamma$ on the $i-$th component
satisfies constraint $\gam_i$ and that
$\gamma(C)$ lies on $k$ general hyperplanes. We use the notation $(\F,\d,k)$ if we don't want
to distinguish the conditions on each compo.nent.

If $\F$ is a closed substack of $\R\R_2(r,d_1,d_2)$ then we denote
$(\F,\gam_1,\gam_2,k,l)$  to be the closure in $\F$
of the locus of maps $\gamma$ such that the
restriction of $\gamma$ on the $i-$th component
satisfies constraint $\gam_i$ and that
$\gamma(C)$ lies on $l$ general hyperplanes, and that
$\gamma(A) = \gamma(B)$ lies on $k$ general hyperplanes.
Similary,
we use the notation $(\F,\d,k,l)$ if we don't want
to distinguish the conditions on each component. 

Note that for maps of reducible source curves, tangency
condition include the case where the image of the node
lies on the tangency hyperplane, as the intersection
multiplicity is $2$ in this case.
\section{Counting one-nodal reducible curves in $\proj$ }
 
In this section we discuss how to count maps with
reducible source curves. 
\begin{proposition}
Let $\F_1$ and $\F_2$ be two
families of stable maps with marked point $C$. Let $\gam_1$
and $\gam_2$ be two constraints. Then we have
  $$\#(\F_1 \times_{ev_{C}} \F_2,\gam_1,\gam_2,k) =  \#( \F_1 , \d_1') \cdot \# (\F_2, \d_2')$$
where $\d_i'$ are determined as follows. Let $e_1$  be 
the dimension of the pushforward under $ev_C$ of $(\F_1,\gam_1)$
into $\proj.$ Let $e_2$ be the dimension of the pushforward under 
$ev_C$ of $(\F_2, \gam_2)$ into $\proj$. Then
$\d_i'$ is obtained from $\gam_i$ by adding
a subspace of codimension $e_i$.
\end{proposition}
\bpf Let $\alpha_i$ be the class of ${ev_C}_*(\F_i,\gam_i)$
in the Chow ring of $\proj$. 
Let $h$ be the class of a subspace of codimension
$k$. Then $\#(\F_1 \times_{ev_{C}} \F_2,\gam_1,\gam_2,k)$
is equal to the intersection product $\alpha_1 \cdot \alpha_2 \cdot h$
which is $\deg(\alpha_1) \cdot \deg(\alpha_2).$ To compute
$\deg(\alpha_i)$, we intersect $\alpha_i$ with a subspace
of codimension $e_i$, thus
$$\deg(\alpha_i) =  \#( \F_i , \d_i')$$
which proves the proposition. \epf \\ \\
The following lemma is useful because it allow us
to express the tangency condition on maps
of reducible curves in terms of tangency conditions
on maps of each component and condition on the node.
\begin{lemma}
 Let $\X_1,\X_2$ be stacks of stable maps into $\proj$.
Assume each map in each family carries
 at least one marked point $C$.
 Let $\X = \X_1 \times_{ev_C} \X_2$
. Let $\T$ be the tangency divisortangenttangent
on $\X$, and $\T_i$ be the pull-back of
the tangency divisor on the $i-$th component.
Then on $\X$ we have this divisorial equation:
$\T = \T_1 + \T_2 + 2\L_C.$
\end{lemma}
\bpf
Let $\C$ be a general curve in $\X$. $\C$ has the following description. There is a family of
nodal curves over $\C,$ $\pi : S \to \C$ such that $S$ is the union of two families of nodal curves $X_1,X_2$ 
along a section $s : \C \to S$. The section $s$ represents the marked point $C$ of each family. There
is also a map $\mu : S \to \proj$
such that the restriction of $\mu$ on each fiber is an element (a map)
of $\X_1 \times_{ev_C} \X_2$. Now choose a general hyperplane $H$ in $\proj.$
Then the restriction of the tangency
divisor $\T$ on $\C$ is the branched divisor of the map $\pi : \mu^{-1}(H)
 = \D \to C$. This map is a $d_1 + d_2$ sheet covering of $\C$. 
The ramification points of this map come from three sources :
\begin{itemize}
\item The ramification points on $\mu^{-1}(H)_{|X_1}.$
\item The ramification points on $\mu^{-1}(H)_{|X_2}.$
\item The intersections $\mu^{-1}(H) \cap s.$
\end{itemize}
The first two sources contribute to the pull backs $\T_1 \cdot \C$ and $\T_2 \cdot \C$ respectively. 
The intersections points $\mu^{-1}(H) \cap s$ correspond 
precisely to the maps $\gamma$ with $\gamma(C) \in H.$ These points
are the nodes of the curve $\D$, because through each of them, there are two branches : one 
from $\mu^{-1}(H)_{|X_1}$, one from $\mu^{-1}(H)_{|X_2}.$ If $P \in \D $ is
 one of such points, then the branched divisor of
$\pi$ contains $\pi(P)$ with multiplicity $2$. Thus we have 
$\T \cdot \C = \T_1 \cdot \C + \T_2 \cdot \C + 2\L_C\cdot \C$.
\epf 

Using the lemma, we can ``expand" the tangency
conditions on $\F_1\times_{ev_C}\F_2$ until we have
tangency conditions only on each individual
component.
\begin{proposition}
 Let $\d$ be a constraint and let $\d_l$
be the constraint obtained from $\d$  by
removing $l$ tangency conditions. Then
we have the following equality :
\begin{eqnarray*}
\# (\F_1 \times_{ev_C} \F_2, \d,k) &=&  \sum_{l=0}^{\d(0)} 2^l \binom{\d(0)}{l} \sum_{\gam_1\gam_2 = \d_l} \# (\F_1\times_{ev_C} \F_2,\gam_1,\gam_2,k+l). \\
\end{eqnarray*}
\end{proposition} 
\bpf There are $(^n_l)$ ways to remove $l$ tangency conditions. 
Doing this results in a codimension $k+l$ condition on the node (the image of $C$)
, and the multiplicity is $2^l.$ \epf 

Applying the proposition to the family $\R\R_2(r,d_1,d_2)$ we have :
\begin{corollary}
$$\#( \R\R_2(r,d_1,d_2),\d,k,k') = \sum_{l=0}^{\d(0)}
2^l\binom{\d(0)}{l}\sum_{\gam_1 \gam_2 = \d_l}\#(\R\R_2(r,d_1,d_2),\gam_1,\gam_2,k,k'+l).$$
\end{corollary}

\section{Counting Rational Space Curves With Special Tangent Conditions}

 In this section, we will describe the algorithm
 counting rational space curves with 
special tangent conditions in $\proj.$ 
 Let $\X = \mbar_{0,\setA}(r,d)$ throughout this section.
 Following the notation in \cite{idq} let
$\H$ be the incident divisor (incident to a codimension 
$2$ subspace), and let $\K^{A,j}$ 
be the boundary divisor of $\mbar_{0,\{A\}}(r,d)$ whose points
represent reducible curves in which the component 
containing $A$ is mapped with degree $j$. 
 The main difficulty when we have multiple special tangent
conditions is excess intersection: any special tangent
divisor $\W_A^M$ passes through the locus of maps 
$\gamma$ where $\gamma(A)$ is not a smooth
point of the image. However, we have the following
result that helps us reduce the number of special
tangent divisors in our computation.
\begin{proposition}
 Any  characteristic number of rational curves with
$l \geq r-1$ special tangent conditions is
expressible in terms of characteristic numbers 
of rational curves with at most $r-2$ special
tangent conditions.
\end{proposition}
 Proof of this statement will be given
in section 5. \epf

Thus, we only need to care about excess intersection
locus in codimension at most $r-2$.  The following
proposition lists all components of this locus.
\begin{proposition}
 Let $S_n$ be the closure of locus of maps $\gamma$ in
$\X$ such that the source curve has $n+1$ components, 
and the component containing $A$, called the principal
component, is incident with $n$ other components.
Moreover, $\gamma$
contracts the principal component. Then
$S_2,\ldots,S_{r-2}$  are the components of codimension
at most $r-2$ of the excess intersection locus
of the special tangent divisors. Furthermore, $S_n$
contributes to the excess intersection only if there are at least
$2n-2$ special tangent conditions. In particular,
only $S_n$'s with $2n \leq r$ are relevant
in counting curves with special tangent conditions.
\end{proposition}
\bpf
 Let $\gamma$ be a map in $\X$ such that
$\gamma(A)$ is not a smooth point of 
its image. If $\gamma$ does not
contract the component of the source curve containing $A$
then $\gamma(A)$ is at least a
nodal singularity. Maps of this type
vary in a family of codimension at
least $r-1$. Thus if $\gamma$
belongs to a component (of the excess
locus) of codimension at most $r-2$,
$\gamma$ must contract the component
of the source curve that containts $A$.
For a multi-index $I(d,n) = (d_1,
\ldots,d_n)$ with $\sum_id_i =d$, let $\V_{I(d,n)}$ be
 $\prod_i \mbar_{0,\setA}(r,d_i)$ where the
product is taken over the evaluation maps
$ev_A$. It is easy to see that each component
of $S_n$ is a finite quotient of a $\mbar_{0,n+1}\times \V_{I(d,n)}$,
where $\mbar_{0,n+1}$ is the moduli
space of genus zero stable curves  with
$n+1$ marked points.
Now $\mbar_{0,n+1}$ is of dimension $n-2$, which means
the "enumerative codimension" of $S_n$ is 
$n-2$ less than its codimension, hence is $2n-2$.
Since we will only need to count
rational curves with at most $r-2$ special tangent
conditions, only $S_n$ in which $2n-2 \leq r-2$,
 or equivalently $2n \leq r$, is relevant.
\epf \\ \\
$$\includegraphics[width = 40mm]{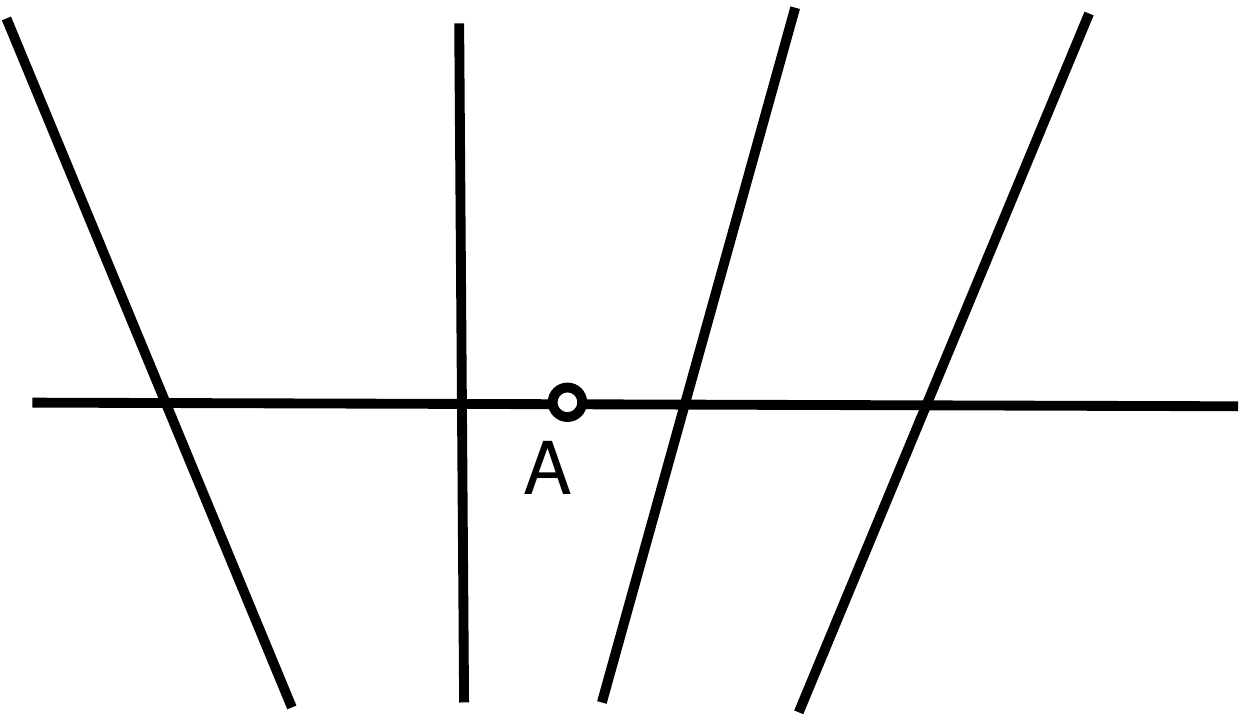}$$ \\
$$\text{Fig 3. ${S_4}$ }$$

 We will blow up $S_n$'s in order to
discount the excess contribution. The above proposition provides 
us with
an useful guideline. In $\mathbb {P}^3$, no blowup is needed. One blowup
of $S_2$ is needed for $\mathbb P^4$ and $\mathbb P^5$. More generally,
we need one more blowup for each increase by two 
in the dimension of the ambient space. In the rest of this section,
we provide explicit formula for the cases $\mathbb P^3,
\mathbb P^4, \mathbb P^5$, which only requires at most
one blowup as expect.  \\ \\
{\bf Case 1:} Counting rational curves with
one special tangent condition in $\proj, r\geq 3.$ \\ \\
We can express the special tangent divisor as 
linear combinations of boundary divisors and
incident divisors, as shown in the following lemma.
\begin{lemma} The following equality holds in the group $A^1(\X) \otimes \rationals$, for $r > 2$ :
$${\mathcal W}_A = 2\L_A + \psi_A$$
where $\psi_A$ is the psi-class. In particular, we have
$$ \W_A = \left(2-\frac{2}{d} \right)\L_A + \frac{1}{d^2}\H + \sum_{j=1}^{j<d}\frac{(d-j)^2}{d^2}\K^{A,j} $$
\end{lemma}
\bpf We use the method as described in \cite{idq}, intersecting the two sides
of the equations with a general curve $\C$ in $\X$. Let $\gamma$ 
denote the image of $\C$ under the evaluation map $ev_A$. Let $M$ be the 
codimension $2$ subspace 
 in $\proj$ corresponding to the special tangent condition $\W_A$.
Beccause $\C$  is a general curve, we can assume
$\gamma$ is smooth.
Let $L$ be a general line in $\proj$,
and let $\pi_M : \proj - M  \to L$
be the projection
 onto $L$ from $M$. Let $\phi_A$ be the line bundle 
on $\gamma $ described as follows. For each point $p \in \gamma$,
$ev_A^{-1}(p)$ is a map $\alpha \in \C$. The fibre of $\phi_A$
over $p$ is then the tangent vector to the image of $\alpha$
at $\alpha(A)$. Let $R$ be the zero scheme of the bundle map
$\phi_A \to  {\pi_M}^*(T_L)$, with $T_L$ being the tangent
bundle of $L$. Geometrically, $R$ represents the 
locus pf points $p \in \gamma$, such that 
the map $ev_A^{-1}(p)$ satisfies special tangent
condition with respect to the subspace $M$. Thus
$$\deg R = R \cap [\gamma] = \C \cap \W_A.$$
We have 
$$\deg R = -c_1(\phi_A) + \deg({\pi_M}_{|\gamma})c_1(T_L).$$
Now $c_1(T_L) = 2$[class of a point], and $\deg({\pi_M}_{|\gamma}) = \deg \gamma
 = \L_A \cap \C$. The pullback of $\phi_A$ by $ev_A$ is isomorphic
to the line bundle on $C$ obtained by attaching to each map
the tangent vector at $A$ to the source curve. Hence
$-c_1(\phi_A)\cap \gamma = -c_1(ev_A^*(\phi_A)) \cap \C = \psi_A \cap  \C$
is the usual psi-class. In short, we have
$${\mathcal W}_A = 2\L_A + \psi_A.$$
The second equality follows from the fact that $\psi_A = - \pi_*(s_A^2)$ on
$\mbar_{0,\seta}(r,d)$ and Lemma 2.2.2 in \cite{idq}. \epf \\ \\
{\bf Case 2:} Counting rational curves with two special tangent 
conditions in $\proj, r\geq 4.$ \\ \\
Let $\pi: \wt{\X} \to \X$
the blowup of $\X$ along $S_2$. Let $S_2^j$ be the component
of $S_2$ with degree partition $(j,0,d-j)$, and let
$E_2^j$ be the corresponding exceptional divisor. We have
that $S_2^j$ is a $\mathbb{Z}_2$-quotient of 
$\R\R(j,d-j)$. A general
element $E_2^j$ has following geometric interpretation: it is a pair
$(\gamma,l)$ where $\gamma$ is a map in $\R\R(j,d-j)$ ,
and $l$ is a line in $\proj$. $l$ must lie on the plane $(l_1,l_2)$
where $l_i$ is the projective tangent line to the image (under $\gamma$)
of the
$i$-th component at the image (under $\gamma$) of $A$ (here we
use $A$ to denote the node of the family $\R\R(j,d-j)$, instead
of using $C$ as in the definition in Section 2.2, but this
does not change anything). For each divisor $\D$
of $\X$, let $\wt{\D}$ be its proper transformation. The next
lemma allows us to compute the class $\pi_*(\wt{\W}_A^2)$
\begin{lemma} 
The following equality holds in $A^2(\X) \otimes \mathbb Q$:
\begin{eqnarray*}  
\pi_* (\wt{\W}_A^2) &=& \left( 2 - \frac{2}{d} \right)\W_A\L_A + \frac{1}{d^2}\W_A\H +
 \sum_{j=1}^{j<d}\frac{(j-d)^2}{d^2}\pi_*(\wt{\W}_A \wt{\K}^{A,j}) + \sum_{j=1}^{j\leq d/2} \frac{2j^2 - 2jd}{d^2}S_2^j
\end{eqnarray*}
The class $\pi_*(\wt{\W}_A \wt{\K}^{A,j})$ is the class of the closure
of the locus of maps with reducible source curves, where
the restriction onto the component containing
$A$ satisfies one special tangent condition.
\end{lemma}

 Counting
maps in $\pi_*(\wt{\W}_A \wt{\K}^{A,j})$ is doable by Lemma $4.1$ and
results in section $3$. Counting maps in $S_2^j$
is equivalent to counting maps in $\R\R(j,d-j)$ which
is also doable by results in section $3$.\\ \\
\bpf We pull back the main equation
of Lemma 4.3:
$$ \pi^*\W_A = \left(2-\frac{2}{d} \right)\wt{\L}_A + \frac{1}{d^2}\wt{\H} + \sum_{j=1}^{j<d}\frac{(d-j)^2}{d^2}\pi^*\K^{A,j} $$
$\pi^*\W_A = \wt{\W}_A + \sum_j E_2^j$ and $\pi^*\K^{A,j} = \wt{\K}^{A,j} + m_jE_2^j$ 
where $m_j$ is $1$ if $j \neq d-j$ and $2$ if $j=d-j$. Rearranging the terms, we have
$$\wt{\W}_A= \left(2-\frac{2}{d} \right)\wt{\L}_A + \frac{1}{d^2}\wt{\H} 
+ \sum_{j=1}^{j<d}\frac{(d-j)^2}{d^2}\wt{\K}^{A,j} +  \sum_{j=1}^{j\leq d/2} \frac{2j^2 - 2jd}{d^2}E_2^j$$
Now it is obvious that $\pi_*(\wt{\W}_AE_2^j)= S_2^j$. Multiply the 
above equation with $\wt{\W}_A$ and pushforward
yields the desired equation. \epf \\ \\
Using Lemma $4.4$, we can reduce
a counting problem involving
two special tangent conditions
into various counting problems
involving at most one special tangent condition. \\ \\
{\bf Case 3:} Counting rational curves
with three special tangent conditions
in $\proj, r\geq 5$. \\ \\
View $\R\R(j,d-j)$ as $\mbar_{0,\setA}(r,j) \times_{ev_A} \mbar_{0,\setA}(r,d-j)$.
Let $\W^{(i)}$ be the pullback of the special
tangent divisor of the $i$-th factor. Let
$p:\R\R(j,d-j) \to S_2^j$ be the natural
projection. We have the following lemma.
\begin{lemma}
 The following equality holds in $A^3(\X)\otimes \rationals$:
\begin{eqnarray*}  
\pi_* (\wt{\W}_A^3) &=& \left( 2 - \frac{2}{d} \right)\pi_*(\wt{\W}^2_A)\L_A + \frac{1}{d^2}\pi_*(\wt{\W}_A)^2\H +
 \sum_{j=1}^{j<d}\frac{(j-d)^2}{d^2}\pi_*(\wt{\W}^2_A \wt{\K}^{A,j}) \\
&+& \sum_{j=1}^{j\leq d/2} \frac{2j^2 - 2jd}{d^2}\pi_*(\wt{\W}_A^2E_2^j)
\end{eqnarray*}
$\pi_*(\wt{\W}^2_A \wt{\K}^{A,j})$ is the closure
in $\X$ of the locus of maps with reducible
source curves, where the restriction of
the map on the component containing $A$
satisfies two special tangent conditions.
Counting maps in this locus is doable
by Lemma $4.4$ and results in section
$3$. Furthermore, for any constraint $\d$ we have
$$(\pi_*(\wt{\W}_A^2E_2^j),\d) = (\W^{(1)} + \W^{(2)},\d)$$
if both sides are finite.
\end{lemma}
\bpf Only the last equality needs proving. Because
the constraint $\d$ cuts out a one-dimensional
family on $\R\R(j,d-j)$, proving the equality is
an intersection theory problem on a $\mathbb P^1$-bunlde 
over a curve. We reformulate the problem
as follows. Let $\F_1$ be a one-dimensional 
family of projective rational curves of degree $j$
with a marked point $A$. We associated with
$\F_1$ the line bundle $l_1$ which is the
line bundle of the projective tangent lines
at $A.$ Similarly, we have $\F_2$ and $l_2$,
where curves in $\F_2$ have degree $d-j.$
Let $\C = \F_1 \times_{ev_A} \F_2$, which is
a curve ($\F_i$'s are choosen so that
$\C$ is not empty). Let $\P$ be the projectivization
of $l_1 \oplus l_2$. Thus $\pi: \P \to \C $
is a rank-one projective bundle. A general
element of $\P$ is a pair of curve-line $(\gamma,l)$
with $\gamma \in \C$ and $l \subset (l_1,l_2).$
Let $\W$
be the divisor on $\P$ define as follows. For
a general codimension $2$ subspace $M\in \proj$,
a pair $(\gamma,l) \in \P$ is in $\W$ if and only 
if $l \subset M$. We have a natural inclusion
$\F_i = P(l_i) \subset \P$, with $P(l_i)$ being the projectivization
of the line bundle $l_i$. Let $\D$ be the canonical
line bundle on $\P$, and let $\G$ be the pullback
of a point $\pi^{(-1)}(p)$ for any $p\in \C$.
 With this reformulation,
the equality that we need to prove becomes
$$\W^2 = \W\F_1 +\W\F_2$$
Let $a_i = - c_1(l_i) \cdot \C$. We have
$$\F_1  = \D + \pi^*(c_1(\phi_2) \cap C) = \D - a_2\G$$
hence 
$$\deg(\F_1^2) =  \deg(\pi_*(\D^2 - 2a_2\D\G + a_2^2\G^2))  = \deg (s_1(F) \cap C) - 2a_2 = a_1 + a_2 - 2a_2 = a_1 - a_2$$
which means that $\F_1^2  = a_1 - a_2$ as $F_1^2$ is of dimension $0$ in the 
Chow ring of $\P$. 
Similarly $\F_2^2 = a_2 - a_1$, thus $F_1^2 + F_2^2 = 0$. Now let $\W =a\F_1 + b\G$. Then we have 
$\W\G = 1 = a(\F\G) \Rightarrow a = 1$.
Now we have $\W\F_1 = \F_1^2 + b \Rightarrow b = \W\F_1 - \F_1^2$. That leads to $\W^2 = 2\W\F_1 - \F_1^2$.
Similarly $\W^2 = 2\W\F_2 - \F_2^2$. Add the two equalities together we have
$$\W^2= \frac{1}{2}(2\W\F_1 + 2\W\F_2 -\F_1^2 - \F_2^2) = \W\F_1 + \W\F_2.$$ \epf \\ \\
Using Lemma $4.5$, we can reduce
a counting problem involving three
special tangent conditions
into various counting problem
involving at most $2$ special
tangent conditions. \\ \\

We end this section with some examples.
\begin{example}
How many  conics in $\mathbb P^3$ passing through $3$ points, that have a marked point $A$
 which must  lie  on a fixed line $M$, and 
that the tangent line at $A$ to the curve
passes through a fixed line $L$? The answer is $1$.
\end{example}
\bpf Because the three points that the conic passes through
determine its plane $H$, this problem reduces to
an enumerative problem in $\mathbb{P}^2$ : how many conics in $\mathbb{P}^2 $ 
that pass through $3$ points and is tangent to a line at a fixed point? The
answer is therefore $1$. Now we will compute this number in a different
way, using Lemma $4.3.$ Let $\d = (0,0,0,3)$, and $\d'=(0,0,1,3)$.
 We need to compute $\# ((\mbar_{0, \{A\}}(3,2), \d),\L^2_AW_A)$.
On $\mbar_{0,\{A\}}(3,2)$, there is one boundary divisor, 
$\K= (\emptyset, 1 \sep \{A \}, 1)$,  which parametrize pair of lines
 intersecting at one point, and
the marked point $A$ is on one of them. Using lemma $4.3$ 
we have
$$\W_A = \L_A + \frac{\H}{4} + \frac{\K}{4}$$
Thus
\begin{eqnarray*}
 \# ((\mbar_{0, \{A\}}(3,2), \d),\L^2_AW_A) &=& \# ((\mbar_{0, \{A\}}(3,2),\d),\L^3_A) + \frac{1}{4}\# 
((\mbar_{0, \{A\}}(3,2), \d'),\L^2_A) \\
                                           && + \frac{1}{4} \# ((\K,\d),\L^2_A)  \\
                                           &=& 0 + \frac{1}{4} + \frac{1}{4}3 = 1.
\end{eqnarray*}
 The first  "$\#$" term of the right hand side  
is the number of conics in $\proj$ 
passing through $4$ points. The second  "$\#$" term is  
the number of conics  in $\proj$ passing through
$3$ points and $2$ lines. The last "$\#$" term
is the number of pair of lines in $\proj$ with one common point,
that pass through $3$ points, and that the component
with the marked point $A$ intersect a line at $A$. \epf
\begin{example}
 There are $2$ conics in $\mathbb{P}^4$ satisfying the following
conditions. The conics pass through $3$ points and a plane, and
there is a marked point $A$ on the curve, the projective
tangent line at which passes through $2$ other planes.
\end{example}
\bpf Again, the three point conditions determine the plane $H$
for the conics. Thus in fact we have a plane curve counting problem.
The conics must pass through $4$ points (the plane condition
now become point condition), and the tangent line at
$A$ must pass through $2$ other points on the plane $H$. Thus
the problem is equivalent to counting plane conics
through $4$ points and tangent to $1$ line, thus the 
answer is two. We must show that
$$\#((\mbar_{0,{A}}(4,2),\d),\W_A^2) = 2 $$
with $\d = (0,1,0,0,3)$. From the proof of 
Lemma $4.4$ we have
$$\wt{\W}_A= \wt{\L}_A + \frac{\wt{\H}}{4} 
+ \frac{\wt{\K}^{A,1}}{4} - \frac{E_2^j}{2}$$
Multiply the equation with $\wt{\W}_A$, pushforward
and integrate against $(\mbar_{0,\{A\}}(4,2),\d)$ we have
\begin{eqnarray*}
 \#((\mbar_{0,{A}}(4,2),\d),\W_A^2) &=& \#((\mbar_{0,{A}}(4,2),\d),\W_A\L_A) + \frac{1}{4}((\mbar_{0,{A}}(4,2),\d'),\W_A) \\
&+& \frac{1}{4}\#((\K^{A,1},\d),\W_A) - \frac{1}{2} \#(E^j_2,\d) \\
 &=& = 3 + \frac{2}{4} + 0 - \frac{3}{2} = 2                                   
\end{eqnarray*}
where $\d' = (0,2,0,0,3)$.	
We list below several numbers of curves with special
tangent conditions in $\mathbb{P}^3, \mathbb{P}^4, \mathbb{P}^5$.
The special class $(a,b)$ means the marked point
as a codimension $a$ condition and there are
$b$ special tangent conditions.
\\ \\
\begin{center}
 \begin{tabular} { | l | c | c | l  |} 
\hline
Degree & Condition & Special Classes & Numbers \\
\hline
Cubic & $(1,2,3)$ & $(3,1)$ & 34 \\
\hline
Cubic & $(4,2,2)$ & $(2,1)$ & 4736 \\
\hline
Quartic & $(7,2,3)$ & $(1,1)$ & 35131904 \\
\hline
Quintic & $(4,4,6)$ & $(0,1)$ & 280111872 \\
\hline
Quintic & $(2,2,7)$ & $(2,1)$ &  352176 \\
\hline
Sextic & $(3,4,7)$ & $(3,1)$ & 340403776 \\
\hline
 \end{tabular}
 $$ \text{Table 1. Some enumerative numbers with special class in $\mathbb{P}^3$} $$
\end{center}

\begin{center}
 \begin{tabular} { | l | c | c | l | } 
\hline
Degree & Condition & Special Classes & Numbers \\
\hline
Conic & $(1,1,2,1)$ & $ (1,2)$ & 38 \\
\hline
Cubic & $(2,1,1,3)$ & $(1,2)$ & $980$ \\
\hline  
Quartic & $(2,2,1,4)$ & $(2,2)$ & $37792$ \\
\hline
Quintic & $(3,3,1,5)$ & $(2,2)$ & $31565232$ \\
\hline
Sextic &$(3,3,4,5)$ &$(1,2)$ & $49679646304$ \\
\hline
 \end{tabular}
$$ \text{Table 2. Some enumerative numbers with special classes in $\mathbb{P}^4$} $$
\end{center}

\begin{center}
 \begin{tabular} { | l | c | c | l | } 
\hline
Degree & Condition & Special Classes & Numbers \\
\hline
Conic & $(1,1,1,0,2)$ & $(0,3)$ & $20$ \\  
\hline
Cubic & $(1,1,1,2,2)$ & $(0,3)$ & $1240 $ \\
\hline
Quartic &$(2,3,1,2,2)$ & $(3,3)$ & $1181400$ \\
\hline
Quintic & $(2,2,3,4,2)$ & $(0,3)$ & $ 1654232816 $ \\
\hline
 \end{tabular}
 $$ \text{Table 3. Some enumerative numbers with special classes in $\mathbb{P}^5$} $$
\end{center}

\section{Counting curves in $\R\R_2(r,d_1,d_2)$}

First we need a result about the Chow ring
of $Bl_{\D}(\proj \times \proj)$, which is
the blowup of $\proj \times \proj$ along
the diagonal. For details of the derivation, we
refer the readers to \cite{dn2}.
\begin{proposition}
The Chow ring of $Bl_{\D}(\proj \times \proj)$ is generated by
$h,k$, the hyperplane class of the first and second factor,
and the exceptional divisor $e$ with the following relations :
\begin{eqnarray*}
h^{r+1} &=& k^{r+1} =0, \\ he&=& ke, \\ e^r &=& \sum_{i>0}^{i < r}(-1)^{i-1}(^{r+1}_{ \ i})h^ie^{r-i} + \sum_{i\geq 0}^{i\leq r}h^ik^{r-i}.
\end{eqnarray*}
\end{proposition} 
\noindent
{\bf Example.}
The following are the third relation in the case $r=1,2,3,4$:
\begin{eqnarray*}
 e &=& h + k. \\
 e^2 &=& 3he - (h^2+hk+k^2). \\
 e^3 &=& 4he^2 - 6h^2e + (h^3 +h^2k + hk^2 + k^3). \\
 e^4 &=& 5he^3 - 10h^2e^2 + 5h^3e - (h^4 + h^3k + h^2k^2 + hk^3 + k^4).
\end{eqnarray*} \epf

 Recall that $\R\R_2(r,d_1,d_2)$ is a substack of
$\mbar_{0, \{A,C\}}(r,d_1) \times_{ev_C} \mbar_{0, \{B,C\}}(r,d_2)$ 
of maps $\gamma $ such that  $\gamma(A) = \gamma(B)$.
 We rephrase the problem of counting maps in $\R\R_2(r,d_1,d_2)$
 as follows : 

{\it Given two families $\F_1$ and $\F_2$ of maps of rational curves with
two marked points $A,C$. How many times
a map $\gamma_1$ from $\F_1$ and a map $\gamma_2 $ from $\F_2$ intersect
in such a way that :  
\begin{itemize}
 \item $\gamma_1(A) = \gamma_2(A)$ and $\gamma_1(C) = \gamma_2(C).$
 \item $\gamma_i(A)$ lies on a fixed linear space of codimension $p$.
 \item $\gamma_i(C)$ lies on a fixed linear space of codimension $q$.
\end{itemize} }

We consider the evaluation map
$$ev_{AC}: \F_i \longrightarrow(\proj\times\proj)$$
Let $T_i$ be the closure in $Bl_{\D}(\proj \times \proj)$
of $ev_{AC}(\F_i)$.
Let $h,k$ be the hyperplane classes of the first
and second factor in $Bl_{\D}(\proj \times \proj).$
Then the answer to our enumerative problem
above is the intersection number
$$T_1T_2h^pk^q$$
where the product is evaluated in the Chow ring
of $Bl_{\D}(\proj \times \proj)$. ($T_i$ parametrizes
ordered pair of points on the curves in $\F_i$. The blowup is
to prevent us from counting in the case where two points
run into each other).

To count maps in $\R\R_2(r,d_1,d_2)$
satisfying the constraint $(\d,p,q),$
 we first consider all the partitions
$\d  = \gam_1\gam_2$, and for each such partition,
assign constraint $\gam_i$ to the $i$-th component.
If $\d(0) \neq 0$, meaning if there are tangency
conditions, we also have to distribute the tangency
conditions over each component first, in the
sense of Proposition $3.3$.
 Then the constraint
$\gam_1$ cuts out a family $\F_1$ on $\mbar_{0,\{A,C\}}(r,d_1)$.
Similarly, $\gamma_2$ cuts out a family $\F_2$ on 
 $\mbar_{0,\{A,C\}}(r,d_2)$. Let
$T_i$ be the closure of $ev_{AC}(\F_i)$ in $Bl_{\D}(\proj \times \proj)$
. We then calculate the product
$$T_1T_2h^pk^q$$
in the Chow ring $A^*(Bl_{\D}(\proj \times \proj))$. Then we take
the sum over all partitions $\d = \gam_1\gam_2$ to
get the number of maps $\# (\R\R_2(r,d_1,d_2), \d ,p,q).$ We need a result to calculate
the classes of $T_i$ in $A^*(Bl_{\D}(\proj \times \proj)).$
The following lemma is useful:
\begin{lemma}
 Let $\F$ be a family of stable maps in $\mbar_{0,\{A,C\}}(r,d)$
such that $A,C$ moves freely, that is, the forgetful map $\mbar_{0,\{A,C\}}(r,d) \to \mbar_{0,0}(r,d)$
has fibre dimension $2$.
Let $T$ be the closure in $Bl_{\D}(\proj \times \proj)$ of the image of $\F$
 under  the evaluation map 
$ev_{AC}: \F \to \proj \times \proj$. Let
$\G$ be the family of stable maps in $\mbar_{0,\{A\}}(r,d)$ 
that is the image of $\F$ under the forgetful morphism
$\mbar_{0,\{A,C\}}(r,d) \to \mbar_{0,\{A\}}(r,d)$. Assume 
$\dim T \leq 2r$. Then we have 
\begin{itemize}
 \item For $m,n$ such that $m+n = \dim T$ :
 $$Th^mk^n= \# (\F,\L_A^m\L_C^n). $$
 \item For $m$ such that $m + 1 = \dim  T$ :
$$Th^me = \# (\G,\L_A^m). $$
 \item For $m,n$ such that $m+n = \dim T$ , we
have $$Th^me(h+k-e)^{n-1}= \#(\G,\L_A^m\W_A^{(n-1)}). $$
\end{itemize}
\end{lemma}
\bpf The first equality is trivial. The number
$Th^mk^n$ is the number of maps  $ \gamma \in \F$ such
that $\gamma(A)$ belongs to $h$ hyperplanes, and
that $\gamma(C)$ belongs to $k$ hyperplanes. That is
precisely the number $\#(\F,\L_A^m\L_C^n).$ The second
equality follows from the fact that multiplying
with $e$ is the same as replacing the family
$\F$ by the family $G$.

Now we
 prove the third equality. Let $$[x_0:x_1:\cdots:x_n] \times [y_0:y_1:\cdots:y_n]$$
be a homogeneous coordinate system of $\proj \times \proj$.
Let $H$ be the hypersurface $$x_0y_n = x_ny_0$$ in $\proj \times \proj.$
$H$ contains $\D$ with multiplicity one and $T = h +k$ in $A^*(\proj \times \proj)$, hence the proper
transformation $\widetilde{H}$ of $H$ in $Bl_{\D}(\proj \times \proj)$
satisfies $$\widetilde{H} = h + k -e.$$ 

 Let us examine what it means
to intersect ${T}$ with $e$ and $\widetilde{H}.$ Let
$\pi : Bl_{\D}(\proj \times \proj) \to \proj \times \proj$
be the blow up, and let $S = \pi(T)$. We have a map $\gamma: S \to S \cap \D$ defined
as folows. For each point
$x \in S$, let $P_x$ be the subspace $ \{p\} \times \proj \subset \proj \times \proj$,
where $\{p\} \in \proj$ is chosen so that $x \in P_x$.
The intersection $S \cap P_x$
is a genus zero curve $f_x$ in $P_x,$ and $\gamma$ maps the entire curve 
$f_x$ onto $x$. The intersection $H \cap P_x$ is a hyperplane in $P_x$
which is the span of $x$ and the codimension $2$ subspace $x_0=y_0 = 0.$
Then for a point $y \in {T}$ with $\pi(y) = x$, we have $ y \in {T} \cap e \cap \wt{H}$
 iff $f_x$, as a curve in the projective space $P_x$
 is tangent to $H_x$ at $x.$
Thus intersecting with $\widetilde{H}$ (after intersecting with $e$) has the
effect of imposing one special tangent condition on the family 
$\G.$ It follows that intersecting with $n-1$ instances of $\wt{T}$ 
has the effect of imposing $n-1$ special tangent conditions.
\epf 
$$\includegraphics[width = 70mm]{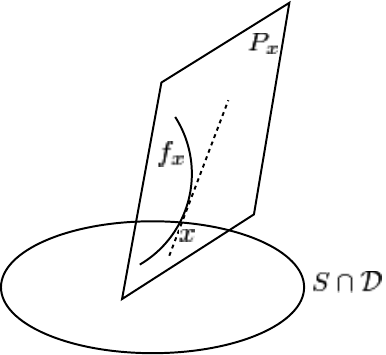}$$
$$\text{Fig 4.}$$\\ \\
Now we have enough to be able compute the class of $T = {ev_{AC}}_*(\F)$
in $A^*(Bl_{\D}(\proj \times \proj)).$ The formal statement of that fact
is the following proposition, whose proof is trivial.
\begin{proposition}
 Let $T \in A^*(Bl_{\D}(\proj \times \proj))$ be a class
of codimension $d, 0 \leq d \leq 2r$ . Then the following intersection products
determine $T$ :
\begin{itemize}
\item $Th^mk^n$ with $0 \leq m \leq r, 0\leq n \leq r.$ 
\item $Th^{m}e(h+k-e)^n$ with $0\leq m \leq r, 0\leq n \leq r-2$.
\item $Th^{d-1}e$.
\end{itemize}
with $m,n$ appropriately choosen so that the 
intersection number is well-defined.
\end{proposition} 
The reason the power $n$ of $h+k-e$ is at most
$r-2$ is because $e^r$ is expressible as polynomials
in $h$ and $k$, so we never need to multiply $T$
with a power of $e$ that is more than $r-1$, in order to determine
$T$.\epf 

In particular, if we know all characteristic numbers
of rational curves with at most $r-2$ special tangent conditions,
then that is enough to count maps in $\R\R_2(d_1,d_2).$ 

{\it Proof of Proposition 4.1.} If the number of special
tangent conditions $l$ is greater than $2r-2$, then the
number is $0$ because the tangent line at $\gamma(A)$
can pass through at most $2r-2$ general codimension
$2$ subspaces. Now assume $l \leq 2r-2$. Let $\d$ be the constraint (beside
the special tangent conditions). Let $\F$ be
$(\mbar_{0,\{A,C\}}(r,d),\d)$ and $T$
be the closure in $Bl_{\D}(\proj \times \proj)$ of 
the image of $\F$ under $ev_{AC}$. We have
$\dim T < 2r.$ If we know
all the characteristic numbers with at most 
$r-2$ special tangent conditions, then
Proposition $5.3$ shows that we can
determine $T$. Then the characteristic
number with constraint $\d$ (and $\L_A^m$) and $l$ 
special tangent conditions is
the intersection number
$Th^me(h+k-e)^l$. \epf
\\ \\
We end the section with some examples.
\begin{example}
 How many pair of lines $(L_1,L_2)$ in $\mathbb P^3$
such that they intersect twice, and that each
of them passes through $3$ lines? The answer is $0.$
\end{example}
 
The answer is obvious because
two distinct lines can never intersect
twice. But our algorithm does not know that.
Let $\d = (0,0,3,0)$. We need to compute
$$\frac{1}{2}\# (\R\R_2(3,1,1), \d,\d). $$
The factor $1/2$ accounts for the fact
that the statement of the problem does
not distinguish the two intersection
points. Let $\F_i$ be the family 
of the lines $L_i$ with a choice of
two marked points $A,C$ on them. 
Let $T_i$ be the pushforward
of $\F_i$ under the evaluation maps
$ev_{AC} : \F_i \to Bl_D(\proj \times \proj).$
$T_1$ is three dimensional, so we can 
assume
$$T_1 = \alpha(h^3+k^3) + \beta(h^2k +hk^2) + \gamma  eh^2 + \mu e^2h. $$
 The coefficients of $h^3$ and $k^3$ must be
the same due to symmetry. Similarly the coefficients
of $h^2k$ and $hk^2$ must be the same.
\begin{eqnarray*}
\alpha &=& \alpha h^3k^3 = T_1k^3 = \# ((\mbar_{0,\{A,C\}}(3,1),\d),\L_A^3) = 0 \\
\beta &=& \beta h^3k^3 = T_1kh^2 = \# ((\mbar_{0,\{A,C\}}(3,1),\d),\L_A^2\L_C) = 2 \\
\mu &=& \mu h^3e^3 = T_1h^2e =  \# ((\mbar_{0,\{A\}}(3,1),\d),\L_A^2) = 2 
\end{eqnarray*}
Computation of $\gamma$ is a little bit lengthier. First we have
\begin{eqnarray*}
\gamma &=& \gamma h^3k^3 = T_1 he^2 - \mu e^4h^2  =  \left (2T_1 h^2e - T_1he(k+k-e) \right) - 4\mu \\
       &=& -2\mu - T_1he(h+k-e). 
\end{eqnarray*}
Now $T_1he(h+k-e) = \# ((\mbar_{0, \{A\}}(3,1),\d),\L_A\W_A)$ is the number of 
lines with a marked point $A$ in $\mathbb P^3$ that pass 
through $3$ lines, such that $A$ lies on
a fixed plane, and such that the tangent line 
at $A$ passes through a general line. This number
is the same as the number of lines passing through
$4$ general lines in $\proj$, which is $2$. Thus
$\gamma = -2\mu - T_1hk(h+k-e) = -4 - 2 = -6$. Therefore
$$T_1 = 2(h^2k+hk^2) - 6h^2e + 2he^2$$
Obviously $T_1 = T_2$, so after a bit of algebra we have
$$T_1T_2 = \left(2(h^2k+hk^2) - 6h^2e + 2he^2 \right)^2= 0.$$ \epf

\begin{example}How many pair of conics-twisted cubics 
in $\mathbb{P}^5$ intersecting at two nodes, with the first node being on a fixed
hyperplane and the second node being on a fixed $3-$space, such that 
the conic passes through one $3-$space, one general
plane, one general line, one general point, and the 
cubic passes through two general $3-$spaces,
one general plane, one general
line, two general points? 
The answer is $956$.
\end{example}

Let $\gam_1 = (0,0,1,1,1,1,0)$
and $\gam_2 = (0,0,2,1,1,2,0)$. We need to compute
$$\#(\R\R_2(5,2,3), \gam_1,\gam_2,1,2).$$
 Let $\F_1$ be a family 
of lines conics in $\mathbb P^5$ with a choice of
two marked points $A,C$ on them, such
that the conics satisfy $\gam_1.$ Let 
$\F_2$ be the a family of twisted cubics
in $\mathbb P^5$ with a choice of two marked 
points $A,C$ on them, such that
the cubics satisfy $\gam_2.$ Let
$T_i$ be the pushforward of $\F_i$
under ${ev_{AC}}_*$ onto
the Chow ring $A^*(Bl_{\D}(\mathbb P^5 \times \mathbb P^5)).$
The we need to compute the intersection product
$hk^2T_1T_2.$
Using Lemma $5.2$ and Proposition $5.3$, we can find the classes
of $T_i$ to be :
\begin{eqnarray*}
 T_1&=& 2h^4 + 6h^3k + 8h^2k^2 + 6hk^3 + 2k^4 - 42h^3e + 29h^2e^2 - 9he^3 + e^4 \\
 T_2&=& 45h^3 + 88h^2k + 88hk^2 + 45k^3 - 308h^2e+ 140he^2 -23e^3
\end{eqnarray*}

Using proposition $5.1$, we can calculate the product:
\begin{eqnarray*}
&&(2h^4 + 6h^3k + 8h^2k^2 + 6hk^3 + 2k^4 - 42h^3e + 29h^2e^2 - 9he^3 + e^4) \\
&\times& ( 45h^3 + 88h^2k + 88hk^2 + 45k^3 - 308h^2e+ 140he^2 -23e^3)hk^2 = 956.
\end{eqnarray*}
\epf 

Some numbers;
\begin{center}
 \begin{tabular} {| l | l | c | c | c | l |}
\hline 
  {\bf Degree} & {\bf Degree} & {\bf Constraint} & {\bf Constraint} & {\bf Nodes} & {\bf Number} \\
\hline
  Conic & Conic & $(2,3,1)$ & $(2,3,1)$ & $(0,0)$ & $ 3360 $ \\
\hline
  Conic & Cubic & $(2,3,1)$ & $(3,4,1)$ & $(1,1)$ & $614656 $ \\
\hline
  Line & Quartic & $(0,1,0)$ &$(3,4,3)$ & $(2,2)$ & $570752$\\
\hline
  Cubic & Cubic & $(3,3,2) $ &$(1,4,2)$ &$(0,3)$ & $963360$\\
\hline
  Conic & Quartic & $(3,3,1) $ & $(0,6,4)$ & $(0,0)$ & $2253312$ \\
\hline
 \end{tabular}
 $$ \text{Table 4. Some enumerative numbers of pair of rational curves in $\mathbb{P}^3$}$$
\end{center}

\begin{center}
 \begin{tabular} {| l | l | c | c | c | l |}
\hline
 {\bf Degree} & {\bf Degree} & {\bf Constraint} & {\bf Constraint} & {\bf Nodes} & {\bf Number} \\
\hline
 Conic & Conic & $(1,1,2,1)$ & $(0,0,0,3)$ & $(0,0)$ & $4$ \\ 
\hline 
 Conic & Cubic & $(1,2,1,1)$ & $(1,1,2,2)$ & $(1,2)$ & $4816$ \\
\hline
 Line & Conic & $(0,1,1,0)$  & $(1,1,1,1)$ & $(1,2)$ & $18$ \\
\hline
 Cubic & Cubic & $(3,1,0,3)$ & $(3,1,0,3)$ & $(1,1)$ & $2297664$ \\
\hline 
 \end{tabular}
 $$ \text{ Table 5. Some enumerative numbers of pair of rational curves in $\mathbb{P}^4$} $$
\end{center}

\begin{center}
 \begin{tabular} {| l | l | c | c | c | l |}
\hline
 {\bf Degree} & {\bf Degree} & {\bf Constraint} & {\bf Constraint} & {\bf Nodes} & {\bf Number} \\
\hline
 Conic & Conic & $(0,0,0,2,1)$ & $(0,0,0,2,1)$ & $(1,1)$ & $2$ \\
\hline
 Conic & Cubic & $(1,0,1,0,2)$ & $(1,1,0,1,3)$ & $(0,0)$ & $144$\\
\hline
 Line & Quartic & $(0,0,0,0,1)$ & $(2,0,0,2,3)$ & $(1,3)$ & $844$ \\
\hline
 Cubic & Cubic & $(3,4,1,1,1)$ & $(2,1,1,2,1)$ & $(1,2)$ & $1027324928$ \\
\hline
 \end{tabular}
 $$ \text{ Table 6. Some enumerative numbers of pair of rational curves in $\mathbb{P}^5$.} $$
\end{center}

\section{Counting rational nodal curves in $\proj$}

First we gave a recursion counting incidence-only
characteristic numbers of rational nodal curves
(with condition on the node) in $\proj$. 
\begin{theorem}
Let
$\d$ be a constraint that $\d(0)= 0.$ Let
$k = \d(r+1).$
Choose a subspace $u$ in $\d$ which is not a hyperplane, such that the dimension of $u$ is largest possible. Then choose any two other subspaces
$s,t$ in $\d$. The following constraints are derived from $\d$ : \\
0) $\wt{\d}$ by removing $u,s,t$ from $\d.$ \\
1) $\d_0$ by replacing $u$ with two subspaces : a hyperplane $p$ and a subspace $q$ such that $p \cap q = u.$ \\
2) $\d_1$ is derived from $\d_0$, by replacing $p$ and $s$ with $p\cap s.$ \\
3) $\d_2$ is derived from $\d_0$, by replacing $q$ and $t$ with $q\cap t.$ \\
4) $\d_3$ is derived from $\d_0$, by replacing $s$ and $t$ with $s\cap t.$ 
\\ \\
If $\gam$ is a set of linear spaces, and $a$ and $b$ are two linear spaces, 
denote $\gam^{(a,b)}$  the  set obtained
from $\gam$ by adding  $a $ and $b.$
 Then the following formula holds :
 \begin{eqnarray*}
\#(\n(r,d),\d) &=& - \sum_{d_1+d_2 = d}^{\gam_1\gam_2 = \wt{\d}} \binom{\wt{\d}}{\gam_1}\#(\n\R(r,d_1,d_2), \gam_1^{(s,t)}, \gam_2^{(p,q)},0)   \\
&-&  \sum_{d_1+d_2=0}^{\gam_1\gam_2 = \wt{\d}} \binom{\wt{\d}}{\gam_1} \#(\n\R(r,d_1,d_2), \gam_1^{(p,q)}, \gam_2^{(s,t)},0) \\
&-&  2\sum_{d_1+d_2 =d}^{\gam_1\gam_2 = \wt{\d}} \binom{\wt{\d}}{\gam_1}\#(\R\R_2(r,d_1,d_2),\gam_1^{(p,q)}, \gam_2^{(s,t)},k,0) \\
&+&  \sum_{d_1+d_2 =d}^{\gam_1\gam_2 = \wt{\d}}\binom{\wt{\d}}{\gam_1} \#(\n\R(r,d_1,d_2),\gam_1^{(q,t)}, \gam_2^{(p,s)},0) \\
&+&  \sum_{d_1+d_2 =d}^{\gam_1\gam_2 = \wt{\d}} \binom{\wt{\d}}{\gam_1}\#(\n\R(r,d_1,d_2),\gam_1^{(p,s)}, \gam_2^{(q,t)},0) \\
&+&  2\sum_{d_1+d_2 =d}^{\gam_1\gam_2 = \wt{\d}} \binom{\wt{\d}}{\gam_1}\#(\R\R_2(r,d_1,d_2),\gam_1^{(p,s)}, \gam_2^{(q,t)},k,0 ) \\
&-& \# (\n(r,d),\d_3) + \# (\n(r,d),\d_1) + \# (\n(r,d),\d_2).
\end{eqnarray*}
Furthermore, $\d_1,\d_2,\d_3$ are all of lower rank than that of $\d.$ 
Here $\binom{\alpha}{\beta} = \prod \binom{\alpha(i)}{\beta(i)}$ for any two tuples $\alpha,\beta$
having the same length.
\end{theorem}
\bpf 
Let $S$ be a set of markings that is in one-to-one
correspondence $\mu : \d_0 \to S$ with the linear spaces in $\d_0$.
Let $\X$  be the moduli space $\mbar_{0, \{ A,B\}
\cup  S}(r,d)$, and let $\n^{(S)}(r,d)$ be the
closure  in $\X$ of the locus of maps $\gamma$
such that $\gamma(A) = \gamma(B).$ Let $\Y$
be the closure in $\n^{(S)}$ of the locus of
maps $\gamma$ such that $\gamma(\mu(m)) \in m$
for all $m \in \d_0.$ Because $\# (\n(r,d),\d)$
is finite, $\Y$ is one-dimensional. We consider
two equivalent divisors on $\X$ :
$$(\{\mu(p),\mu(q) \} \sep \{ \mu(s), \mu(t)\})= (\{ \mu(p), \mu(s) \}\sep \{\mu(q), \mu(t)\}).$$
Let $\K_1 = (\{\mu(p),\mu(q)\} \sep \{\mu(s), \mu(t)\}), $ and let $ \K_2 = (\{\mu(p),\mu(s)\} \sep \{\mu(q), \mu(t)\})$. Then we have
$$\# \left(\Y \cap \K_1\right) = \# \left( \Y \cap \K_2 \right).$$
Let us analyze the left-hand side of the equation. Let
$\gamma$ be a general point of $\Y \cap \K_1$.
Then $\gamma$ is a stable map whose source curve
has two components $C_1,C_2$ joined at a node,
such that $\mu(p), \mu(q) \in C_1$ and
$\mu(s), \mu(t) \in C_2.$
There are several cases to consider:
\begin{itemize}
 \item $\deg \gamma_{|C_1} =0.$ If only $A$ 
or $C$ is on
$C_1$ then by dimension couting we have that
this case has no contribution. If both $A,C$ are
on $C_1$ then the image curve has a cusp, on which
we impose condition like those we impose on $p,q$.
By dimension count again, we also have that
the case has no contribution. The quick reason
is that if a map contracted a component containing
at least $4$ special points (marked or nodes),
then the dimension of the family of image curves
is less than the dimension of the family of maps,
therefore is enumeratively irrelevant.
Now if
$A,B \in C_2$,  $\gamma_{|C_2}$ is a rational
nodal curve and satisfies the constraint $\d$
(but these conditions are marked). The contribution
to $\# (\Y \cap \K_1)$ in this case is $\#(\n(r,d),\d).$
 \item $\deg \gamma_{|C_2} = 0.$ Arguing similarly,
we have that the contribution to $\# (\Y \cap \K_1)$
is $\# (\n(r,d), \d_3)$
 \item $\gamma$ has positive degree $d_i$ component $C_i.$
There are three subcases :
\begin{itemize}
 \item $A,B \in C_1:$ In this case,
$\gamma_{|C_1}$ is a rational nodal curve
and $\gamma_{|C_2}$ is a rational curve. The
contribution in this case is 
 $$\sum_{d_1+d_2=0}^{\gam_1\gam_2 = \wt{\d}}\#(\n\R(r,d_1,d_2), \gam_1^{(p,q)}, \gam_2^{(s,t)},0).$$
\item $A,B \in C_2$ : The contribution is 
$$ \sum_{d_1+d_2 = d}^{\gam_1\gam_2 = \wt{\d}} \#(\n\R(r,d_1,d_2), \gam_1^{(s,t)}, \gam_2^{(p,q)},0). $$
\item $A \in C_1, B\in C_2$ or vice versa.
In this case the image of $\gamma$ is
a curve having two components that intersect
twice at distinguished points. The contribution
is therefore
$$ 2\sum_{d_1+d_2 =d}^{\gam_1\gam_2 = \wt{\d}} \#(\R\R_2(r,d_1,d_2),\gam_1^{(p,q)}, \gam_2^{(s,t)},k,0 ).$$
\end{itemize}
\end{itemize}
We can analyze $\Y  \cap  \K_2$ in the  same
 way and after rearranging the terms, 
we derive the equation in the statement of the 
theorem. \epf

It is now possible to use the results so far to compute the characteristic
number of rational nodal curves.
\begin{theorem}
Let $\d$ be a constraint such that $\d(0)  >0$. Let
$\d(r+1) = k$
Let $\d''$ be the constraint obtained from $\d$
by removing a tangency hyperplane. Let $\d'$ be the constraint
obtained from $\d''$ by adding an incident codimension $2$
subspace. Then
we have the following equality, provided that the
left hand side is finite.

\begin{eqnarray*}
 \#(\n(r,d), \d) &=& \frac{d-1}{d}  \# (\n(r,d), \d') \\ 
&+& \sum_{d_1+d_2 = d} \big( \#  (\n\R(r,d_1,d_2),\d'') + \#( \R\R_2(r,d_1,d_2),\d'',k,0) \big).
\end{eqnarray*}
\end{theorem}
{\bf Warning :} if $\d(0) \neq 0$ then those summands above involving
reducible curves contain (twice) the case where the node
is mapped to a tangency hyperplane. Also, in computing
those summands, one needs to consider all possible splitting
of constraints over two components (see Proposition $3,3$ and
Corollary $3,4$). 

\bpf We have the following equality of divisors on 
$\mbar_{0,\{A.B\}}(r,d)$
$$\T = \frac{d-1}{d}\H + \sum_{d>0}^{j \leq d/2}\frac{j(d-j)}{d}(j,d-j).$$
For a proof of this see \cite{idq}, Lemma $2.3.1$.
Thus
\begin{eqnarray*}
\# (\n(r,d), \d) &=& \#  \left ((\n(r,d),\d''),\T \right ) \\
               &=& \frac{d-1}{d} \# \left ( (\n(r,d),\d''),\H \right) + \sum_{j>0}^{j \leq d/2} \# \left( \n(r,d) \cap (j,d-j), \d'' \right).
\end{eqnarray*}
Now we will analyze $\#( \n(r,d) \cap (j,d-j), \d'')$.
A general map $\gamma \in \n(r,d) \cap (j,d-j)$ has
two-component source curve.
There are two cases:

\begin{itemize}
 \item $A,B$ belong to a same component. The contribution 
is $\#(\n\R(j,d-j),\d'') + \#(\n\R(d-j,j),\d'')$ if $j < d-j$ depending
on whether $A,B$ are in the component of lower or
higher degree. If $j = d-j,$ the contribution is
just $\#(\n\R(j,d-j),\d'').$
\item $A,B$ belong to different components. The contribution
is $2\#( \R\R_2(j,d-j),\d'',k,0)$ if $j < d-j$ and is
$\#( \R\R_2(j,d-j), \d'',k,0)$ if $j =d-j.$
\end{itemize}
Sum up all possibilities, we derive the formula
in Theorem $6.2$. \epf \\ \\

Calculation of $\#(\R\R_2(r,d_1,d_2),\d'',k,0)$ should
make use of Corollary $3.4.$
One point worth mentioning when counting
rational nodal curves with tangency conditions
and with condition on the node is that
maps with degree $2$ do contribute enumeratively.
Rational nodal curves with degree two are
rational degree two covers of $\mathbb P^1$ with
a marked point specified as the node. For these
maps, having a hyperplane passing through the 
branched points count as tangency. 

From characteristic number of rational nodal curves, it is easy
to get characteristic number of rational nodal curves. Let $m = \d(0)$,
 and $\d_i$ be the constraint received by removing $i$ tangency
conditions and replace them by a codimension $i$ on the node. Then
we have the number of elliptic curves with fixed $j-$ invariant,
with $j$ generic, of degree $d$ in $\proj$ satisfying constraint 
$\d$ denoted $\#(\J(r,d),\d)$, is :
$$\#(\J(r,d),\d) = \sum_{i=0}^m 2^i \binom{n}{i}\#(\n(r,d),\d_i).$$

Now we give several numerical examples. We recover
all previously known numbers in literature.
The characteristic numbers of plane nodal cubics were
computed in \cite{luf}. The charactersitic numbers of elliptic plane
curves with fixed $j-$ invariant were computed in \cite{char}. Charactersitic numbers of rational plane cubics
in $\mathbb{P}^3$ were computed in \cite{hmx}. Let $N,N_l,N_p$ be 
the family of rational nodal curves, rational nodal curves
with the node on a fixed line, rational nodal curves
with the node on a fixed point. Similarly, we denote
$N_{s},N_{b}, N_{f}$ for the same family with the node on a fixed plane,
 a fixed $3-$space, or a fixed $4-$space. The following tables list the characteristic
numbers of such families and of elliptic curves with fixed $j-$
invariant (denoted by $\J$). Below are tables of characteristic numbers of such families of
low degree ($2,3,4,5$). In some tables, we put some point conditions
so that the numbers are small enouch to fit in the table. The only other
conditions are tangency, and top incident condition. For example, in 
the table for quartics in $\mathbb{P}^4$, the curves must pass through
$2$ points, the other conditions are combination of tangency and
incident to planes.

\begin{center}
 \begin{tabular} { | c | l | l | l | l |}
\hline
  $\#$ tang & $N$ & $N_l$ & $N_p$ & $\J$  \\
\hline
 $0$& $0$ & $0$  & $0$ &  $0$ \\
\hline
 $1$ & $0$ & $0$  & $0$ &  $0$ \\ 
\hline
 $2$ & $0$ & $2$  & $1$ &  $0$ \\
\hline
 $3$ &  $0$ & $3$  & $3/2$ & $12$ \\
\hline
 $4$ & $0$ & $3/2$ &      & $48$ \\
\hline
 $5$ & $0$ &       &      & $75$ \\
\hline
 \end{tabular}
$$ \text{Table 7. Plane conics.} $$
\end{center}

\begin{center}
 \begin{tabular} {| c | l | l| l| l |}
 \hline
  $\#$ tang & $N$ & $N_l$ & $N_p$ & $\J$  \\
 \hline
  $0$ & $12$  &   $6 $  & $1$ & $12$ \\
 \hline
  $1$  &    $36$ &  $22$ & $4$ & $48$ \\
 \hline
  $2$  &  $100$ &  $80$  &  $16$  & $192$ \\
 \hline
  $3$  &   $240$ &  $240$  & $52$  & $768$ \\
 \hline
  $4$  &   $480$ &  $604$  & $142$ & $2784$ \\
 \hline
  $5$  &  $712$ &   $1046$ & $256$ & $8832$ \\
 \hline
  $6$  &   $756$ &  $1212$ & $304$ & $21828$ \\
 \hline
  $7$  &   $600$ &  $1000$ &  & $39072$ \\
 \hline
  $8$  &   $400$ &     &  & $50448$\\
\hline
 \end{tabular}
 $$ \text{Table 8. Plane cubics.} $$
\end{center}

\begin{center}
 \begin{tabular} {| c | l | l| l| l |}
 \hline
  $\#$ tang & $N$ & $N_l$ & $N_p$ & $\J$  \\
 \hline
  $ 0 $ &  $ 1860 $ &  $ 768 $ &  $ 96 $ &  $ 1860 $  \\ 
\hline
$ 1 $ &  $ 6552 $ &  $ 2952 $ &  $ 384 $ &  $ 8088 $  \\ 
\hline
$ 2 $ &  $ 21600 $ &  $ 10712 $ &  $ 1448 $ &  $ 33792 $  \\ 
\hline
$ 3 $ &  $ 65328 $ &  $ 35616 $ &  $ 4992 $ &  $ 134208 $  \\ 
\hline
$ 4 $ &  $ 178272 $ &  $ 106752 $ &  $ 15516 $ &  $ 497952 $  \\ 
\hline
$ 5 $ &  $ 429120 $ &  $ 281348 $ &  $ 42416 $ &  $ 1696320 $  \\ 
\hline
$ 6 $ &  $ 886632 $ &  $ 633972 $ &  $ 99024 $ &  $ 5193768 $  \\ 
\hline
$ 7 $ &  $ 1515960 $ &  $ 1166352 $ &  $ 187248 $ &  $ 13954512 $  \\ 
\hline
$ 8 $ &  $ 2097648 $ &  $ 1705856 $ &  $ 279152 $ &  $ 31849968 $  \\ 
\hline
$ 9 $ &  $ 2350752 $ &  $ 1986672 $ &  $ 329496 $ &  $ 60019872 $  \\ 
\hline
$ 10 $ &  $ 2184480 $ &  $ 1893528 $ &  & $ 92165280 $  \\ 
\hline
$ 11 $ &  $ 1745712 $ &  &  &$ 115892448 $  \\ 
\hline
 \end{tabular}
$$ \text{Table 9. Plane quartics.} $$
\end{center}

\begin{center}
  \begin{tabular} {| c | l | l| l| l | l |}
 \hline
  $\#$ tang & $N$ & $N_s$ & $N_l$ & $N_p$ & $\J$  \\
\hline
   $ 0 $ &  $ 0 $ &  $ 0 $ &  $ 0 $ &  $ 0$ &  $ 0 $  \\ 
\hline
$ 1 $ &  $ 0 $ &  $ 0 $ &  $ 0 $ &  $ 0 $ &  $ 0 $  \\ 
\hline
$ 2 $ &  $ 0 $ &  $ 16 $ &  $ 8 $ &  $ 2 $ &  $ 0 $  \\ 
\hline
$ 3 $ &  $ 0 $ &  $ 24 $ &  $ 12 $ &  $ 3 $ &  $ 96 $  \\ 
\hline
$ 4 $ &  $ 0 $ &  $ 20 $ &  $ 10 $ &  $ 7/2 $ &  $ 384 $  \\ 
\hline
$ 5 $ &  $ 0 $ &  $ 10 $ &  $ 5 $ &   &  $ 840 $  \\ 
\hline
$ 6 $ &  $ 0 $ &  $ 5 $ &   &   &  $ 1200 $  \\ 
\hline
$ 7 $ &  $ 0 $ &   &   &   &  $ 1470 $  \\ 
\hline
  \end{tabular}

$$ \text{Table 10. Conics in $\mathbb{P}^3$.} $$
\end{center}
\begin{center}
  \begin{tabular} {| c | l | l| l| l | l |}
 \hline
  $\#$ tang & $N$ & $N_s$ & $N_l$ & $N_p$ & $\J$  \\
\hline
  $ 0 $ &  $ 12960 $ &  $ 5040 $ &  $ 904 $ &  $ 72 $ &  $ 12960 $  \\ 
\hline
$ 1 $ &  $ 29520 $ &  $ 13120 $ &  $ 2512 $ &  $ 216 $ &  $ 39600 $  \\ 
\hline
$ 2 $ &  $ 61120 $ &  $ 32048 $ &  $ 6568 $ &  $ 612 $ &  $ 117216 $  \\ 
\hline
$ 3 $ &  $ 109632 $ &  $ 64608 $ &  $ 13904 $ &  $ 1384 $ &  $ 332640 $  \\ 
\hline
$ 4 $ &  $ 167616 $ &  $ 107072 $ &  $ 23904 $ &  $ 2524 $ &  $ 849024 $  \\ 
\hline
$ 5 $ &  $ 214400 $ &  $ 144960 $ &  $ 33304 $ &  $ 3732 $ &  $ 1890240 $  \\ 
\hline
$ 6 $ &  $ 230240 $ &  $ 162760 $ &  $ 38432 $ &  $ 4656 $ &  $ 3625440 $  \\ 
\hline
$ 7 $ &  $ 211200 $ &  $ 155288 $ &  $ 37808 $ &  $ 5112 $ &  $ 5994096 $  \\ 
\hline
$ 8 $ &  $ 170192 $ &  $ 130048 $ &  $ 32864 $ &  $ 5424 $ &  $ 8631120 $  \\ 
\hline
$ 9 $ &  $ 124176 $ &  $ 98352 $ &  $ 25664 $ &   &  $ 11038224 $  \\ 
\hline
$ 10 $ &  $ 85440 $ &  $ 70880 $ &   &   &  $ 12875520 $  \\ 
\hline
$ 11 $ &  $ 56960 $ &   &   &   &  $ 14422080 $  \\ 
\hline

 \end{tabular}
  $$ \text{Table 11. Cubics in $\mathbb{P}^3$.} $$
\end{center}

\begin{center}
  \begin{tabular} {| c | l | l| l| l | l |}
 \hline
  $\#$ tang & $N$ & $N_s$ & $N_l$ & $N_p$ & $\J$  \\
\hline

$ 0 $ &  $ 247191840 $ &  $ 61582704 $ &  $ 7487280 $ &  $ 402216 $ &  $ 247191840 $  \\ 
\hline
$ 1 $ &  $ 519424512 $ &  $ 138566640 $ &  $ 17469840 $ &  $ 975192 $ &  $ 642589920 $  \\ 
\hline
$ 2 $ &  $ 1034619648 $ &  $ 295896480 $ &  $ 38636160 $ &  $ 2242512 $ &  $ 1618835328 $  \\ 
\hline
$ 3 $ &  $ 1932171072 $ &  $ 588656160 $ &  $ 79348512 $ &  $ 4785408 $ &  $ 3920405760 $  \\ 
\hline
$ 4 $ &  $ 3353134848 $ &  $ 1079389056 $ &  $ 149728320 $ &  $ 9378160 $ &  $ 9020858112 $  \\ 
\hline
$ 5 $ &  $ 5361957120 $ &  $ 1808973504 $ &  $ 257515200 $ &  $ 16752296 $ &  $ 19509189120 $  \\ 
\hline
$ 6 $ &  $ 7841572992 $ &  $ 2752793920 $ &  $ 401264800 $ &  $ 27140752 $ &  $ 39298619520 $  \\ 
\hline
$ 7 $ &  $ 10431095808 $ &  $ 3788712880 $ &  $ 564734880 $ &  $ 39830752 $ &  $ 73227372288 $  \\ 
\hline
$ 8 $ &  $ 12599060192 $ &  $ 4716456320 $ &  $ 718744512 $ &  $ 53161088 $ &  $ 125665152480 $  \\ 
\hline
$ 9 $ &  $ 13851211968 $ &  $ 5333385216 $ &  $ 831757440 $ &  $ 65099040 $ &  $ 198307833792 $  \\ 
\hline
$ 10 $ &  $ 13948252800 $ &  $ 5522229504 $ &  $ 883153920 $ &  $ 74131776 $ &  $ 288227491200 $  \\ 
\hline
$ 11 $ &  $ 12986719872 $ &  $ 5292561600 $ &  $ 870495360 $ &  $ 79929312 $ &  $ 387635041920 $  \\ 
\hline
$ 12 $ &  $ 11309818368 $ &  $ 4757882880 $ &  $ 807883200 $ &  $ 84550992 $ &  $ 486058242048 $  \\ 
\hline
$ 13 $ &  $ 9330496512 $ &  $ 4070594880 $ &  $ 715629312 $ &   &  $ 574243507200 $  \\ 
\hline
$ 14 $ &  $ 7394421888 $ &  $ 3381893376 $ &   &   &  $ 648194719872 $  \\ 
\hline
$ 15 $ &  $ 5703866880 $ &   &   &   &  $ 715490590080 $  \\ 
\hline
 \end{tabular}
 $$ \text{Table 12. Quartics in $\mathbb{P}^3$.}$$
\end{center}

\begin{center}
  \begin{tabular}{| c | l | l| l| l | l |}
 \hline
  $\#$ tang & $N$ & $N_s$ & $N_l$ & $N_p$ & $\J$  \\
\hline
   $ 0 $ &  $ 2987074368 $ &  $ 597069288 $ &  $ 59293632 $ &  $ 2757288 $ &  $ 2987074368 $  \\ 
\hline
$ 1 $ &  $ 6654861504 $ &  $ 1393675584 $ &  $ 142403568 $ &  $ 6890568 $ &  $ 7849000080 $  \\ 
\hline
$ 2 $ &  $ 14302171008 $ &  $ 3141287760 $ &  $ 330349200 $ &  $ 16691344 $ &  $ 20114047872 $  \\ 
\hline
$ 3 $ &  $ 29534616768 $ &  $ 6800411520 $ &  $ 736077600 $ &  $ 38978688 $ &  $ 50113244448 $  \\ 
\hline
$ 4 $ &  $ 58394890752 $ &  $ 14081928256 $ &  $ 1569037056 $ &  $ 87466348 $ &  $ 120947061888 $  \\ 
\hline
$ 5 $ &  $ 110164217088 $ &  $ 27795971008 $ &  $ 3189343752 $ &  $ 188200508 $ &  $ 281761911168 $  \\ 
\hline
$ 6 $ &  $ 197654921184 $ &  $ 52144209544 $ &  $ 6165495488 $ &  $ 387843208 $ &  $ 631585386720 $  \\ 
\hline
$ 7 $ &  $ 336286484448 $ &  $ 92755042440 $ &  $ 11312688400 $ &  $ 765476504 $ &  $ 1358700870672 $  \\ 
\hline
$ 8 $ &  $ 541376364848 $ &  $ 156271230640 $ &  $ 19684719200 $ &  $ 1449944208 $ &  $ 2800306366128 $  \\ 
\hline
$ 9 $ &  $ 823917940992 $ &  $ 249556959696 $ &  $ 32520764016 $ &  $ 2653490208 $ &  $ 5526457857888 $  \\ 
\hline
$ 10 $ &  $ 1186459103808 $ &  $ 379132252128 $ &  $ 51221741472 $ &  $ 4769939328 $ &  $ 10455705197568 $  \\ 
\hline
$ 11 $ &  $ 1621483284864 $ &  $ 552185368704 $ &  $ 77488852608 $ &   &  $ 19030887269760 $  \\ 
\hline
$ 12 $ &  $ 2114474172288 $ &  $ 783085854720 $ &   &   &  $ 33559605535872 $  \\ 
\hline
$ 13 $ &  $ 2648546358528 $ &   &   &   &  $ 58098921777408 $  \\ 
\hline
  \end{tabular}
 $$ \text{ Table 13. Quintics in $\mathbb{P}^3$, passing through $3$ points.} $$
\end{center}

\begin{center}
 \begin{tabular} {|  c | l | l | l | l| l | l |} 
 \hline $\#$ tang & $N$ & $N_b$ & $N_s$ & $N_l$ & $N_p$ & $\J$  \\
 
 \hline
  $ 0 $ &  $ 7833840 $ &  $ 2565720 $ &  $ 468935 $ &  $ 52140 $ &  $ 2865 $ &  $ 7833840 $  \\ 
\hline
$ 1 $ &  $ 14708400 $ &  $ 5294270 $ &  $ 1017980 $ &  $ 119400 $ &  $ 6984 $ &  $ 19839840 $  \\ 
\hline
$ 2 $ &  $ 25085900 $ &  $ 10073080 $ &  $ 2038520 $ &  $ 252192 $ &  $ 15720 $ &  $ 48138720 $  \\ 
\hline
$ 3 $ &  $ 37705920 $ &  $ 16296840 $ &  $ 3416336 $ &  $ 440272 $ &  $ 28924 $ &  $ 110777280 $  \\ 
\hline
$ 4 $ &  $ 49732080 $ &  $ 22491008 $ &  $ 4833312 $ &  $ 644504 $ &  $ 44470 $ &  $ 232897920 $  \\ 
\hline
$ 5 $ &  $ 57643520 $ &  $ 26854560 $ &  $ 5889580 $ &  $ 812540 $ &  $ 59250 $ &  $ 439941120 $  \\ 
\hline
$ 6 $ &  $ 59232320 $ &  $ 28240140 $ &  $ 6319450 $ &  $ 906690 $ &  $ 70854 $ &  $ 745702080 $  \\ 
\hline
$ 7 $ &  $ 54660200 $ &  $ 26636130 $ &  $ 6095150 $ &  $ 916962 $ &  $ 78360 $ &  $ 1141405440 $  \\ 
\hline
$ 8 $ &  $ 45993500 $ &  $ 22938610 $ &  $ 5383586 $ &  $ 858012 $ &  $ 82584 $ &  $ 1593774300 $  \\ 
\hline
$ 9 $ &  $ 35861700 $ &  $ 18337518 $ &  $ 4423952 $ &  $ 755184 $ &  $ 85440 $ &  $ 2055201960 $  \\ 
\hline
$ 10 $ &  $ 26323500 $ &  $ 13808900 $ &  $ 3420200 $ &  $ 626640 $ &  $ 87360 $ &  $ 2480472300 $  \\ 
\hline
$ 11 $ &  $ 18497240 $ &  $ 9949360 $ &  $ 2513120 $ &  $ 480480 $ &   &  $ 2841879120 $  \\ 
\hline
$ 12 $ &  $ 12649200 $ &  $ 6978480 $ &  $ 1786880 $ &   &   &  $ 3137555760 $  \\ 
\hline
$ 13 $ &  $ 8510880 $ &  $ 4808480 $ &   &   &   &  $ 3385230720 $  \\ 
\hline
$ 14 $ &  $ 5673920 $ &   &   &   &   &  $ 3589051200 $  \\ 
\hline
 \end{tabular}
$$ \text{Table 14. Cubics in $\mathbb{P}^4$.} $$
\end{center}

\begin{center}
 \begin{tabular} {|  c | l | l | l | l| l | l |} 
 \hline $\#$ tang & $N$ & $N_b$ & $N_s$ & $N_l$ & $N_p$ & $\J$  \\
 
 \hline
  $ 0 $ &  $ 264271032 $ &  $ 61079694 $ &  $ 8388348 $ &  $ 749421 $ &  $ 34860 $ &  $ 264271032 $  \\ 
\hline
$ 1 $ &  $ 493716948 $ &  $ 120918936 $ &  $ 17290038 $ &  $ 1630488 $ &  $ 81252 $ &  $ 615876336 $  \\ 
\hline
$ 2 $ &  $ 878434848 $ &  $ 228232116 $ &  $ 33980664 $ &  $ 3390452 $ &  $ 181836 $ &  $ 1395663984 $  \\ 
\hline
$ 3 $ &  $ 1479817080 $ &  $ 405964896 $ &  $ 62797160 $ &  $ 6629800 $ &  $ 383672 $ &  $ 3062685600 $  \\ 
\hline
$ 4 $ &  $ 2353692768 $ &  $ 678089744 $ &  $ 108738088 $ &  $ 12151512 $ &  $ 761888 $ &  $ 6469681248 $  \\ 
\hline
$ 5 $ &  $ 3530480992 $ &  $ 1063566824 $ &  $ 176508768 $ &  $ 20905076 $ &  $ 1429930 $ &  $ 13101001152 $  \\ 
\hline
$ 6 $ &  $ 4995675728 $ &  $ 1569827616 $ &  $ 269290448 $ &  $ 33879818 $ &  $ 2556172 $ &  $ 25387171536 $  \\ 
\hline
$ 7 $ &  $ 6680908448 $ &  $ 2189197336 $ &  $ 387775734 $ &  $ 51989792 $ &  $ 4399696 $ &  $ 47102511264 $  \\ 
\hline
$ 8 $ &  $ 8472417440 $ &  $ 2900923506 $ &  $ 529920660 $ &  $ 75922720 $ &  $ 7378752 $ &  $ 83878893600 $  \\ 
\hline
$ 9 $ &  $ 10234272948 $ &  $ 3679075344 $ &  $ 691414728 $ &  $ 105627552 $ &  $ 12126048 $ &  $ 143940578328 $  \\ 
\hline
$ 10 $ &  $ 11836475952 $ &  $ 4504817304 $ &  $ 867212688 $ &  $ 138946656 $ &   &  $ 239302639872 $  \\ 
\hline
$ 11 $ &  $ 13167563808 $ &  $ 5374257696 $ &  $ 1054871808 $ &   &   &  $ 387833169936 $  \\ 
\hline
$ 12 $ &  $ 14112721248 $ &  $ 6278297856 $ &   &   &   &  $ 616383262944 $  \\ 
\hline
$ 13 $ &  $ 14531107200 $ &   &   &   &   &  $ 963518793600 $  \\ 
\hline
 \end{tabular}
$$ \text{Table 15. Quartics in $\mathbb{P}^4$ passing through $2$ points.} $$
\end{center}

\begin{center}
 \begin{tabular}{|  c | l | l | l | l| l | l |} 
 \hline $\#$ tang & $N$ & $N_b$ & $N_s$ & $N_l$ & $N_p$ & $\J$  \\
 \hline
 $ 0 $ &  $ 5264130996 $ &  $ 960390870 $ &  $ 105886953 $ &  $ 7801695 $ &  $ 311311 $ &  $ 5264130996 $  \\ 
\hline
$ 1 $ &  $ 10335707556 $ &  $ 1973618742 $ &  $ 224710598 $ &  $ 17371678 $ &  $ 742316 $ &  $ 12256489296 $  \\ 
\hline
$ 2 $ &  $ 19791788388 $ &  $ 3960252460 $ &  $ 465840460 $ &  $ 37911496 $ &  $ 1746624 $ &  $ 28109811168 $  \\ 
\hline
$ 3 $ &  $ 36896035320 $ &  $ 7737537944 $ &  $ 940326944 $ &  $ 80796848 $ &  $ 4041128 $ &  $ 63416490816 $  \\ 
\hline
$ 4 $ &  $ 66880583024 $ &  $ 14699954352 $ &  $ 1845469104 $ &  $ 167905648 $ &  $ 9189708 $ &  $ 140521932288 $  \\ 
\hline
$ 5 $ &  $ 117792292576 $ &  $ 27145486560 $ &  $ 3519654728 $ &  $ 340028520 $ &  $ 20558296 $ &  $ 305497218816 $  \\ 
\hline
$ 6 $ &  $ 201506364736 $ &  $ 48745168872 $ &  $ 6523861268 $ &  $ 670681448 $ &  $ 45308086 $ &  $ 651327035136 $  \\ 
\hline
$ 7 $ &  $ 334871977648 $ &  $ 85223104580 $ &  $ 11759484440 $ &  $ 1287078386 $ &  $ 98524384 $ &  $ 1362231952128 $  \\ 
\hline
$ 8 $ &  $ 540951986840 $ &  $ 145379939744 $ &  $ 20637848154 $ &  $ 2397410108 $ &  $ 211715288 $ &  $ 2797819372056 $  \\ 
\hline
$ 9 $ &  $ 850242885024 $ &  $ 242702404542 $ &  $ 35332114224 $ &  $ 4312424928 $ &   &  $ 5652591017568 $  \\ 
\hline
$ 10 $ &  $ 1301286873156 $ &  $ 397849014300 $ &  $ 59181220928 $ &   &   &  $ 11257978051236 $  \\ 
\hline
$ 11 $ &  $ 1938666465816 $ &  $ 641728301752 $ &   &   &   &  $ 22149199999776 $  \\ 
\hline
$ 12 $ &  $ 2804649121008 $ &   &   &   &   &  $ 43096623642288 $  \\ 
\hline
 \end{tabular}
$$ \text{Table 16. Quintics in $\mathbb{P}^4$ passing through $4$ points.} $$
\end{center}

\begin{center}
 \begin{tabular}{|  c | l | l |  l | l | l| l | l |} 
 \hline $\#$ tang & $N$ & $N_f$ & $N_b$ & $N_s$ & $N_l$ & $N_p$ & $\J$  \\
 \hline
  $ 0 $ &  $ 3580435656 $ &  $ 1034759292 $ &  $ 189136374 $ &  $ 24039939 $ &  $ 2009982 $ &  $ 85745 $ &  $ 3580435656 $  \\ 
\hline
$ 1 $ &  $ 5820250128 $ &  $ 1803057816 $ &  $ 343203840 $ &  $ 45424176 $ &  $ 3974516 $ &  $ 178640 $ &  $ 7889768712 $  \\ 
\hline
$ 2 $ &  $ 8641680264 $ &  $ 2888520852 $ &  $ 572163144 $ &  $ 78755588 $ &  $ 7205344 $ &  $ 341240 $ &  $ 16610457024 $  \\ 
\hline
$ 3 $ &  $ 11507535984 $ &  $ 4048138080 $ &  $ 824350976 $ &  $ 116897472 $ &  $ 11089152 $ &  $ 549128 $ &  $ 33149426688 $  \\ 
\hline
$ 4 $ &  $ 13759570272 $ &  $ 4992894416 $ &  $ 1036797728 $ &  $ 150683904 $ &  $ 14773856 $ &  $ 764324 $ &  $ 61362323712 $  \\ 
\hline
$ 5 $ &  $ 14867247680 $ &  $ 5502189760 $ &  $ 1161050240 $ &  $ 172833416 $ &  $ 17554792 $ &  $ 954832 $ &  $ 104391383040 $  \\ 
\hline
$ 6 $ &  $ 14650427520 $ &  $ 5502894720 $ &  $ 1179603568 $ &  $ 180279708 $ &  $ 19079772 $ &  $ 1102606 $ &  $ 163351745280 $  \\ 
\hline
$ 7 $ &  $ 13303631040 $ &  $ 5066847184 $ &  $ 1104900496 $ &  $ 174051444 $ &  $ 19343536 $ &  $ 1204100 $ &  $ 236503108800 $  \\ 
\hline
$ 8 $ &  $ 11252393152 $ &  $ 4350397184 $ &  $ 967029476 $ &  $ 157723006 $ &  $ 18576208 $ &  $ 1267280 $ &  $ 319397674176 $  \\ 
\hline
$ 9 $ &  $ 8959119120 $ &  $ 3522421644 $ &  $ 799569876 $ &  $ 135605388 $ &  $ 17095224 $ &  $ 1305896 $ &  $ 405992118672 $  \\ 
\hline
$ 10 $ &  $ 6782773704 $ &  $ 2715749316 $ &  $ 629998440 $ &  $ 111418656 $ &  $ 15173120 $ &  $ 1331840 $ &  $ 490193697672 $  \\ 
\hline
$ 11 $ &  $ 4929887760 $ &  $ 2011043040 $ &  $ 476256768 $ &  $ 87775688 $ &  $ 12973792 $ &  $ 1349216 $ &  $ 567210910536 $  \\ 
\hline
$ 12 $ &  $ 3472645440 $ &  $ 1442366496 $ &  $ 347592224 $ &  $ 66354624 $ &  $ 10586880 $ &  $ 1360832 $ &  $ 634363027200 $  \\ 
\hline
$ 13 $ &  $ 2392303152 $ &  $ 1010425424 $ &  $ 246674816 $ &  $ 48224736 $ &  $ 8073728 $ &   &  $ 691172850672 $  \\ 
\hline
$ 14 $ &  $ 1624181888 $ &  $ 696607744 $ &  $ 171675392 $ &  $ 34118336 $ &   &   &  $ 738716078016 $  \\ 
\hline
$ 15 $ &  $ 1092498624 $ &  $ 474968256 $ &  $ 117859840 $ &   &   &   &  $ 778457098944 $  \\ 
\hline
$ 16 $ &  $ 730705920 $ &  $ 321392512 $ &   &   &   &   &  $ 811258656768 $  \\ 
\hline
$ 17 $ &  $ 487137280 $ &   &   &   &   &   &  $ 838048055040 $  \\ 
\hline
 \end{tabular}
$$ \text {Table 17. Cubics in $\mathbb{P}^5$} $$
\end{center}

\begin{center}
 \begin{tabular}{|  c | l | l |  l | l | l| l | l |} 
 \hline $\#$ tang & $N$ & $N_f$ & $N_b$ & $N_s$ & $N_l$ & $N_p$ & $\J$  \\
 \hline
$ 0 $ &  $ 17793468 $ &  $ 4315338 $ &  $ 675729 $ &  $ 82815 $ &  $ 7629 $ &  $ 408 $ &  $ 17793468 $  \\ 
\hline
$ 1 $ &  $ 33892524 $ &  $ 8728578 $ &  $ 1428506 $ &  $ 187086 $ &  $ 18804 $ &  $ 1122 $ &  $ 42523200 $  \\ 
\hline
$ 2 $ &  $ 61915284 $ &  $ 16962956 $ &  $ 2898296 $ &  $ 406116 $ &  $ 44736 $ &  $ 3012 $ &  $ 99532512 $  \\ 
\hline
$ 3 $ &  $ 108109320 $ &  $ 31398264 $ &  $ 5580216 $ &  $ 834384 $ &  $ 100788 $ &  $ 7728 $ &  $ 227691648 $  \\ 
\hline
$ 4 $ &  $ 180450912 $ &  $ 55359984 $ &  $ 10188624 $ &  $ 1618620 $ &  $ 214248 $ &  $ 18948 $ &  $ 507304944 $  \\ 
\hline
$ 5 $ &  $ 288477120 $ &  $ 93327232 $ &  $ 17697268 $ &  $ 2968056 $ &  $ 429304 $ &  $ 44638 $ &  $ 1099292256 $  \\ 
\hline
$ 6 $ &  $ 442955328 $ &  $ 151262244 $ &  $ 29385528 $ &  $ 5155156 $ &  $ 807974 $ &  $ 101692 $ &  $ 2318653056 $  \\ 
\hline
$ 7 $ &  $ 655304328 $ &  $ 237174048 $ &  $ 46930448 $ &  $ 8512992 $ &  $ 1413096 $ &   &  $ 4771225200 $  \\ 
\hline
$ 8 $ &  $ 936129552 $ &  $ 361876128 $ &  $ 72589134 $ &  $ 13497600 $ &   &   &  $ 9605588880 $  \\ 
\hline
$ 9 $ &  $ 1291589856 $ &  $ 539604810 $ &  $ 109323720 $ &   &   &   &  $ 18969484704 $  \\ 
\hline
$ 10 $ &  $ 1716845652 $ &  $ 788940756 $ &   &   &   &   &  $ 36822211764 $  \\ 
\hline
$ 11 $ &  $ 2184938712 $ &   &   &   &   &   &  $ 70374247152 $  \\ 
\hline

 \end{tabular}
$$ \text{Table 18. Quartics in $\mathbb{P}^5$ passing through $3$ points.} $$
\end{center}

\end{document}